\documentclass[12pt,leqno]{article}
\usepackage{amsfonts}
\usepackage{amsmath, amsthm, amsfonts, amssymb, color}
\usepackage{mathrsfs}
\usepackage{color}
\theoremstyle{definition}

\newcommand{\scr}[1]{\mathscr #1}
\definecolor{wco}{rgb}{0.5,0.2,0.3}

\numberwithin{equation}{section} \theoremstyle{remark}

\newcommand{\ua}{\uparrow}

\title{{\bf Estimates for   Invariant Probability Measures of Degenerate SPDEs with Singular and  Path-Dependent Drifts}\footnote{Supported in
 part by  NNSFC (11771326, 11431014, 11726627).} }
\author{
{\bf     Feng-Yu Wang  }\\
\footnotesize{Center of Applied Mathematics, Tianjin University, Tianjin 300072, China}\\
\footnotesize{  wangfy@tju.edu.cn }}
\begin{document}
\allowdisplaybreaks
\def\R{\mathbb R}  \def\ff{\frac} \def\ss{\sqrt} \def\B{\mathbf
B} \def\W{\mathbb W}
\def\N{\mathbb N} \def\kk{\kappa} \def\m{{\bf m}}
\def\ee{\varepsilon}\def\ddd{D^*}
\def\dd{\delta} \def\DD{\Delta} \def\vv{\varepsilon} \def\rr{\rho}
\def\<{\langle} \def\>{\rangle} \def\GG{\Gamma} \def\gg{\gamma}
  \def\nn{\nabla} \def\pp{\partial} \def\E{\mathbb E}
\def\d{\text{\rm{d}}} \def\bb{\beta} \def\aa{\alpha} \def\D{\scr D}
  \def\si{\sigma} \def\ess{\text{\rm{ess}}}
\def\beg{\begin} \def\beq{\begin{equation}}  \def\F{\scr F}
\def\Ric{\text{\rm{Ric}}} \def\Hess{\text{\rm{Hess}}}
\def\e{\text{\rm{e}}} \def\ua{\underline a} \def\OO{\Omega}  \def\oo{\omega}
 \def\tt{\tilde} \def\Ric{\text{\rm{Ric}}}
\def\cut{\text{\rm{cut}}} \def\P{\mathbb P} \def\ifn{I_n(f^{\bigotimes n})}
\def\C{\scr C}      \def\aaa{\mathbf{r}}     \def\r{r}
\def\gap{\text{\rm{gap}}} \def\prr{\pi_{{\bf m},\varrho}}  \def\r{\mathbf r}
\def\Z{\mathbb Z} \def\vrr{\varrho} \def\ll{\lambda}
\def\L{\scr L}\def\Tt{\tt} \def\TT{\tt}\def\II{\mathbb I}
\def\i{{\rm in}}\def\Sect{{\rm Sect}}  \def\H{\mathbb H}
\def\M{\scr M}\def\Q{\mathbb Q} \def\texto{\text{o}} \def\LL{\Lambda}
\def\Rank{{\rm Rank}} \def\B{\scr B} \def\i{{\rm i}} \def\HR{\hat{\R}^d}
\def\to{\rightarrow}\def\l{\ell}\def\iint{\int}
\def\EE{\scr E}\def\Cut{{\rm Cut}}\def\V{\mathbb V}
\def\A{\scr A} \def\Lip{{\rm Lip}}
\def\BB{\scr B}\def\Ent{{\rm Ent}}\def\L{\scr L}\def\div{{\rm div}}
\maketitle

\begin{abstract}   In terms of a nice reference   probability measure,  integrability conditions on the  path-dependent drift are presented for (infinite-dimensional) degenerate   PDEs to have regular positive solutions.  To this end, the corresponding stochastic (partial) differential equations are proved  to possess the weak existence and  uniqueness  of   solutions,  as well as the  existence, uniqueness  and entropy estimates  of  invariant probability measures.   When the reference measure satisfies the log-Sobolev inequality,  Sobolev estimates are derived for the density of invariant probability measures. Some results are new even for non-degenerate SDEs with path-independent drifts. The main results are applied to   nonlinear functional  SPDEs and degenerate functional  SDEs/SPDEs.

\end{abstract} \noindent
 \vskip 2cm

\section{Introduction}

It is well known that hypoelliptic  differential operators with smooth coefficients share similar  properties  with the elliptic ones. For instances, the H\"ormander theorem \cite{H} ensures the    smoothness of heat kernels (see Malliavin \cite{Ma} for a probabilistic proof),  the Index theorem has been proved  by Bismut \cite{B}. See also  \cite{NN, SK1, SK2} and references within for Harnack inequality for hypoelliptic equations and asymptotics of heat kernels.  In this paper, we investigate estimates of positive solutions to (infinite-dimensional) hypoelliptic equations with   singular and path-dependent drifts.

Consider, for instance,  the following second order differential operator in $\R^d$:
$$\scr L := \sum_{i=1}^m X_i^2 +X_0,$$ where $X_0,\cdots, X_m$  are locally bounded  vector fields. A function $\rr\in L_{loc} (\R^d)$   is called a weak solution to $\scr L^*\rr=0$ if
$$\int_{\R^d} (\rr \scr L f)(x)\d x=0,\ \ f\in C_0^\infty(\R^d).$$
Similarly, a locally finite signed measure $\nu$ is called a solution to the equation $\scr L^*\nu=0$ if
$$\int_{\R^d}  \scr L f \d \nu =0,\ \ f\in C_0^\infty(\R^d).$$

When $X_0,\cdots, X_m$ are smooth such that Lie$\{X_1,\cdots, X_m\}$ has rank $d$ (i.e. $\{X_1,\cdots, X_m\}$ satisfies the H\"ormander condition), a positive solution $\rr$ to the equation $\scr L^*\rr=0$ is locally H\"older continuous with respect to the intrinsic distance and satisfies the local Harnack inequality, see \cite{NN} and references within.

When the operator is non-degenerate, i.e. rank$\{X_1,\cdots, X_m\}=d,$ the drift $X_0$ is allowed to be very singular. More precisely, rewrite $\scr L= {\rm tr}(a\nn^2) +b\cdot\nn$, where $b\in L^1_{loc}(\R^d\to \R^d)$ and $a=\si\si^*$  for some $\si\in L_{loc}^1(\R^d\to \R^d\otimes\R^d)$ such that $a\ge \ll I$ for some $\ll\in C(\R^d; (0,\infty))$. If $a$ is differentiable in the distribution sense and $|\nn a|+ |b|\in L_{loc}^p(\R^d)$ for some $p>d$, then
any positive solution to the equation $\scr L^*\nu=0$ for measures has a strictly positive density $\rr \in W^{1,p}_{loc}(\R^d)$, see \cite{BKR, BR} and references within for more results in the literature.
Recently, explicit integrability conditions are presented in \cite{W16} to imply the existence, uniqueness and global regularity estimates on probability solutions to the equation $\scr L^*\nu=0$.

Here, we investigate probability solutions to $\scr L^*\nu=0$ for more general and more singular situations. For instance,   we   consider   differential operators  on the path space $\C:= C([-\tau,0];\R^d)$ for some $\tau>0$ (in the following, $\R^d$ will be extended to an Hilbert space). For any $h,g\in L^2([-\tau,0];\R^d)$, let $\<h,g\>_2=\int_{-\tau}^0\<h(\theta),g(\theta)\>\d\theta.$    We  introduce the class $\F C$ of   cylindrical functions of type
$$\xi\mapsto f\big(\xi(0), \<h_1, \xi\>_2,\cdots, \<h_n, \xi\>_2\big), \ \ n\ge 1, h_i\in C^1([-\tau,0];\R^d).$$ Consider the following  path-dependent operator $\scr L_{a,b}$:   for the above type function $f\in \F C$ and $\xi\in \C$,
\beg{align*} &\scr L_{a,b} f(\xi) :=  \Big(\sum_{i,j=1}^d  a_{ij}(\xi(0))  \pp_i\pp_j    +\sum_{i=1} b_i(\xi) \pp_i\Big) f\big(\cdot, \<h_1,\xi\>_2,\cdots, \<h_n, \xi\>_2\big)(\xi(0))\\
&\quad + \sum_{k=1}^n \big\{\<h_k(0),\xi(0)\>-\<h_k(-\tau),\xi(-\tau)\> -\<h_k', \xi\>_2\big\}\big\{\pp_k f(\xi(0),\cdot) \big\}(\<h_1,\xi\>_2,\cdots, \<h_1,\xi\>_2),\end{align*}
where $a:=(a_{ij})_{1\le i,j\le d}\ge 0$  (maybe degenerate) is $C^2$ but $b:=(b_1,\cdots, b_d): \C\to \R^d$ only satisfies an integrability condition with respect to a probability measure $\mu$. So,  $b$ might be only  $\mu$-a.e.  defined without any continuity. We will construct the Markov semigroup   generated by $\scr L_{a,b}$ (i.e. weak solutions to the corresponding SDE),  and investigate the invariant probability measures. In general, an invariant probability measure $\nu$ of the semigroup solves the equation $\scr L_{a,b}^* \nu=0$ in the sense that
$$ \int_\C (\scr  L_{a,b} f)(\xi)\nu(\d\xi)=0,\ \ f\in \F C.$$

To explain how far we will go beyond the existing study, let us briefly recall the main result in \cite{W16}.
 Consider the following SDE    on $\R^d$:
\beq\label{SDE}\d X(t) =\{Z_0+\si Z\}(X(t))\d t + \si(X(t))\d W(t),\end{equation}  where   $W(t)$ is the $d$-dimensional Brownian motion; $Z: \R^d\to\R^d$ is measurable;   $\si\in C^2(\R^d\to\R^d\otimes\R^d)$ such that $\si(x)$ is invertible for every $x\in\R^d$ and the intrinsic metric
$$\rr_\si(x,y):=\sup\big\{|f(x)-f(y)|:\ f\in C^\infty(\R^d), |\si^*\nn f|\le 1\big\},\ \ x,y\in \R^d$$ is complete; and
\beq\label{Z0}Z_0= \ff 1 2 \sum_{i,j=1}^d \{\pp_j (\si\si^*)_{ij}   - (\si\si^*)_{ij} \pp_j V\}e_i\end{equation} for some $V\in C^2(\R^d)$ and  the standard orthonormal  basis $\{e_i\}_{i=1}^d$ of $\R^d$. Let $\mu_0(\d x) =\e^{-V(x)}\d x$, and let
  $H^{1,2}_\si(\mu_0)$ be  the closure of $C_0^\infty(\R^d)$ under the norm
  $$\|f\|_{H^{1,2}_\si(\mu_0)}:=\big\{\mu_0(|f|^2+|\si^*\nn f|^2)\big\}^{\ff 1 2}. $$
 The  following result is taken from \cite[Theorem 2.1 and Theorem 2.3]{W16}  where   the constant $\kk$   is replaced by $2 \kk$ as the noise therein is $\ss 2 W(t)$ rather than $W(t).$ 

 \beg{thm}[\cite{W16}]\label{T1.1} Let $\si\in C^2(\R^d\to\R^d\otimes\R^d)$ such that $\si(x)$ is invertible for every $x\in\R^d$,   let $Z_0$ be in $\eqref{Z0}$ for some $V\in C^2(\R^d)$, and let $\mu_0(\d x)=\e^{-V(x)}\d x.$ \beg{enumerate} \item[$(1)$] If   for some constant $\vv\in (0,1)$ $$ \int_{\R^d}  \e^{\vv |Z(x)|^2 -\vv^{-1}\rr_\si(0,x)^2  }\mu_0(\d x)<\infty, $$ then for any initial points the SDE $\eqref{SDE}$  has a unique  non-explosive solution, and the associated Markov semigroup  $P_t^Z$  is strong Feller with at most one invariant probability measure.
 \item[$(2)$] Let $\mu_0$ be a probability measure satisfying the (defective) log-Sobolev inequality
\beq\label{LS} \mu_0(f^2\log f^2)\le \kk \mu_0(|\si^*\nn f|^2)+\bb,\ \ f\in C_0^\infty(\R^d), \mu_0(f^2)=1 \end{equation}
   for some   constants $\kk>0,\bb\ge 0.$ If
\beq\label{SIG}\mu_0\big(\e^{\vv \|\si\|^2}+ \e^{\ll |Z|^2}\big):=\int_{\R^d} \big\{\e^{\vv\|\si||^2}+\e^{\ll |Z|^2}\big\}\d\mu_0<\infty\end{equation} holds for some constants
$\vv>0$ and  $\ll> \kk,$ then   $P_t^Z$ has a unique  invariant probability measure $\nu$, which is absolutely continuous with respect to $\mu_0$  and   $\rr:=\ff{\d\nu}{\d\mu_0} $ has a continuous, strictly positive version such that
  $\ss\rr,\,\log\rr \in H^{1,2}_\si(\mu_0)$ with
\beg{equation}\beg{split} \label{SI2}&\mu_0\big(|\si^*\nn\ss{\rr}|^2\big) \le \ff{1}{\ll- \kk} \big\{\log\mu_0(\e^{\ll|Z|^2}) +\bb\big\}< \infty,\\
&\mu_0(|\si^*\nn\log\rr|^2) \le 4 \mu_0(|Z|^2) <\infty.\end{split}\end{equation} \end{enumerate} \end{thm}

Since an invariant probability measure $\nu$ of $P_t^Z$ solves the   equation
$\scr L_Z^*\nu=0,$  due to the integration by  parts formula, Theorem \ref{T1.1}(2) provides regularity estimates on positive  solutions to the singular elliptic PDE
$$ \Big(\ff 1 2 {\rm Tr}(\si\si^* \nn^2)  +(Z_0-\si Z)\cdot\nn\Big) \rr =0.$$

We will  improve and extend the above  assertions in the following four aspects:
\beg{enumerate} \item[$\bullet$] The noise may be degenerate: $|\si^* \nn f|=0$ does not imply $\nn f=0$,  so that the log-Sobolev inequality \eqref{LS} is invalid. Moreover, $\si$ is not necessarily $C^2$-smooth.
 \item[$\bullet$] The drift $Z$ may be path-dependent, for which the corresponding SDE is called functional SDE.
 \item[$\bullet$] The state space may be infinite-dimensional  such that the study   applies to   nonlinear or semilinear functional   SPDEs.
 \item[$\bullet$] Derive  stronger  estimates on the density of  the invariant probability measure. \end{enumerate}
We now introduce the framework of the present study in details.

\paragraph{Reference SDE.} Let $\H$ and $\tt\H$ be two separable Hilbert spaces, and let $\L(\tt \H;\H)$ be the class of bounded linear operators from $\tt\H$ to $\H$. The cylindrical Brownian motion on $\tt\H$ is formally defined by
$$W(t)= \sum_{i\ge 1} B_i(t)\tt e_i,$$ where $\{\tt e_i\}_{i\ge 1}$ is an  orthonormal  basis of $\tt \H$,  and  $\{B_i(t)\}_{i\ge 1}$ are  independent one-dimensional Brownian motions on a complete filtration probability space
$(\OO, \{\F_t\}_{t\ge 0}, \P)$. Let $\V$ be a Banach space densely embedded into $\H$, and let $\V^*$ be its dual space with respect to $\H$. We call
  $\V\subset \H\subset \V^*$   the Gelfand triple.
Consider the following reference SDE   on  $\H$:
 \beq\label{E0} \d X(t)=  Z_0(X(t)) \d t + \si(X(t))\d W(t),\end{equation} where
   $Z_0:\H\to \V^*$ and $\si: \H\to \L(\tt\H;\H)$ are measurable. A continuous adapted process   $X(t)$  on $\H$ is called a (variational) solution to \eqref{E0} with initial value $X(0)$, if
   $$  \E \int_0^t \big\{|\,_{\V^*}\<Z_0(X(s)), v\>_\V|+  |\{\si(X(s))\}^*v|^2\big\}\d s <\infty$$ for all $t\in (0,\infty)$ and $v\in \V$, and   $\P$-a.s.
   $$\<X(t),v\>_\H= \<X(0),v\>_\H+ \int_0^t \,_{\V^*}\<Z_0(X(s)), v\>_\V\,\d s +\int_0^t \big\<\{\si(X(s))\}^*v, \d W(s)\big\>_{\tt\H},\ \ v\in\V, t\ge 0.$$
 See \cite{LR} and reference within for the existence and uniqueness of variational solutions under  framework of monotone SPDEs due to \cite{Par, KZ}. When the initial value $X(0)=x$, we denote the solution by $X^x(t)$. When  the initial value $X(0)$ has distribution $\nu_0$ on $\H$, we also denote the solution by $X^{\nu_0}(t)$ to emphasize the initial distribution.

 When   $\H=\V=\R^d$, the variational solution reduces to   the usual strong (i.e. pathwise) solution of SDEs.    When $\V=\H$ and $Z_0(x)= A x + \tt Z_0(x)$, where   $\tt Z_0\in C(\H;\H)$ and $A$ is a   self-adjoint operator  on $\H$  generating a $C_0$-contraction semigroup $T_t$   such that
 $$\E\int_0^t \Big(\sum_{i\ge 1} |T_{t-s}\si(X(s))\tt e_i|^2\d s+ |T_{t-s}\tt Z_0(X(s))|\Big)\d s<\infty,\ \ t\ge 0,$$
 it coincides with the mild solution  in the sense  of \cite{DZ,DZ2}.

 Throughout the paper, we assume:

 \beg{enumerate}  \item[{\bf (A)}] \emph{Let $Z_0:\H\to \V^*$ and $\si: \H\to \L(\tt\H;\H)$ be measurable such that for any $\F_0$-measurable random variable   $X(0)$ on $ \H$,   the SDE \eqref{E0} has
  a unique variational solution,     and  the associated Markov $P_t^0$ given by
  $$P_t^0f(x):= \E \big[f(X^x(t))\big],\ \  t\ge 0, x\in\H, f\in \B_b(\H)$$ has a unique invariant probability measure $\mu_0$.}\end{enumerate}

Under assumption {\bf (A)}, for any probability measure $\nu_0$ on $\H$, we have
$$\E [f(X^{\nu_0}(t))] =\int_{\H} \E [f(X^x(t))] \nu_0(\d x)= \nu_0(P_t^0 f).$$ In particular,  $\E [f(X^{\mu_0}(t))]=\mu_0(f)$ for $t\ge 0$ and $f\in \B_b(\H)$.

Let  $\tau>0$. For any $\xi\in C([-\tau,\infty); \H)$ and $t\ge 0$, define $\xi_t\in \C:= C([-\tau,0];\H)$ by
$$\xi_t(\theta):= \xi(t+\theta),\ \ \ \theta\in [-\tau, 0].$$ We call $(\xi_t)_{t\ge 0}$ the segment of $(\xi(t))_{t\ge -\tau}$.       For an $\F_0$-measurable random variable $\xi$ on $\C$, let
$$X^\xi(t) = \xi(t)1_{[-\tau,0]}(t) + X^{\xi(0)}(t) 1_{(0,\infty)}(t),\ \ t\ge -\tau,$$ recall that $X^{\xi(0)}(t)$ is the solution to \eqref{E0} with initial value $\xi(0)$.  Let $(X_t^\xi)_{t\ge 0} $ be   the segment process of $(X^\xi(t))_{t\ge -\tau}$, i.e.
\beq\label{D*E}X_t^\xi(\theta)= X^\xi(t+\theta),\ \ \theta\in [-\tau,0].\end{equation}   When  $\xi$ has distribution $\nu$, we also denote $X_t^\xi$ by $X_t^\nu$ to emphasize the initial distribution. Then
\beq\label{34}S_t^0 f(\xi):= \E[f(X_t^\xi)],\ \ \xi\in \C, f\in \B_b(\C), t\ge 0  \end{equation} gives rise to a Markov semigroup $S_t^0$ on $\B_b(\C)$.

 \paragraph{Reference probability measure $\mu$.}  Let $\mu$ be  the distribution of the $\C$-valued random variable  $X^{\mu_0}_\tau$ defined by $X_\tau^{\mu_0}(\theta):=X^{\mu_0}(\tau+\theta), \theta\in [-\tau,0]$. Recall that $X^{\mu_0}(t)$ is the solution to \eqref{E0} with initial distribution $\mu_0$. It is easy to see that $\mu$ is the unique invariant probability measure of the Markov semigroup $S_t^0$.

 Since $X_t^\xi(0):=X^\xi(t)= X^{\xi(0)}(t)$  for   $t\ge 0$,
$$(S_{t+\tau} f)(x):= (S_{t+\tau}^0 f)(\xi),\ \ \xi\in\C, \xi(0)=x, f\in \B_b(\C)$$ provides a family of contractive  linear operators
$(S_{t+\tau})_{t\ge 0}$ from $\B_b(\C)$ to $  \B_b(\H).$    By the Markov property, this implies
\beq\label{3*}\beg{split} &(S_{t+\tau}f)(\xi(0))= (S_{t+\tau}^0 f)(\xi) = \E[(S_\tau^0f )(X_t^\xi)]
 = \E[(S_\tau f)(X_t^\xi(0))]\\
 & =\E[(S_\tau f)(X^{\xi(0)}(t))]= P_t^0(S_\tau f)(\xi(0)),   \ \ \xi\in \C, t\ge 0, f\in \B_b(\C).\end{split}\end{equation}

 \paragraph{Singular and path-dependent  SDE.} Consider the following SDE on $\H$:
 \beq\label{E1} \d X(t) =\big\{Z_0(X(t))  +\si(X(t))Z(X_t)\big\}\,\d t + \si(X(t))\d W(t),\end{equation}
where $Z: \C \to  \tt\H$  and $(X_t)_{t\ge 0}$ is the segment process of   $(X(t))_{t\ge -\tau}.$  Even in the path-independent case, when $\si$ is degenerate or $\H$ is infinite-dimensional, to ensure the strong existence and uniqueness one needs certain continuity conditions on the drift,  see \cite{DFRP, DFRV, W15, WZ, ZN} and references within for details. So, to   investigate   \eqref{E1} by using
 integrability conditions of $Z$ with respect to the reference   measure $\mu$, we only look at the weak solution.

\beg{defn} Let $\xi\in \C$. A $\C$-valued continuous process $(\tt X_t^\xi)_{t\ge 0}$ under a complete filtration probability space
$(\tt\OO, \{\tt\F_t\}_{t\ge 0},\tt \P)$ is called a weak solution of \eqref{E1} starting at $\xi$, if it is $\tt \F_t$-adapted with $\tt X_0^\xi=\xi$ and for some $\tt \H$-cylindrical Brownian motion $\tt W(t)$ on the same probability space  $\tt \P$-a.s.
\beg{align*} \<\tt X^\xi(t),v\>_\H= &\<\xi(0),v\>_\H +\int_0^t \,_{\V^*}\big\<  Z_0(\tt X^\xi(s))+\si(\tt X^\xi(s))Z(\tt X^\xi_s), v\big\>_\V \d s\\
 &+ \int_0^t \big\<\{\si(\tt X^\xi(s))\}^*v, \d \tt W(s)\big\>_{\tt\H},\ \ t\ge 0,v\in\V. \end{align*}   The equation is said to have weak uniqueness if any two weak solutions with same initial point are equal in law. \end{defn}

When \eqref{E1} has weak existence and uniqueness,  let $\P^\xi$ denote the distribution of the weak solution starting at $\xi$, and define
$$(S_t^Zf)(\xi) :=\int_{C([0,\infty);\C)}  f(\eta(t))\ \P^\xi(\d\eta),\ \ f\in \B_b(\C),\xi\in\C, t\ge 0.$$ A probability measure $\nu$ on $\C$ is called an invariant probability measure   of $S_t^Z$,   if $\nu(S_t^Z f)=\nu(f)$ holds for all $t>0$ and $f\in \B_b(\C).$

\

The remainder of the paper is organized as follows. 
Under condition $\mu(\e^{\vv|Z|^2})<\infty$ for some $\vv>0$, we prove  the weak existence and uniqueness of solutions (Section 2)  as well as  the uniqueness of  invariant probability measure (Section 3).
Moreover, the existence of invariant probability measures  and  entropy estimate of the density are   proved in Section 4 using the hyperboundedness of $P_t^0$. Finally, the existence of invariant probability measures and  Sobolev estimates on the  density are addressed  in Section 5 by using the  log-Sobolev inequality \eqref{LS}, for which the H\"ormander condition is adopted. The main results are applied to concrete  models of degenerate functional SDEs/SPDEs. We emphasize that  some estimates  in Section 5 are new even for non-degenerate SDEs with path-independent drifts, see Theorem \ref{T3.1} and Theorem \ref{T3.2} below for details.

\section{Weak solutions}

 Let $\mu(\e^{\vv|Z|^2})<\infty$ for some constant $\vv>0$. We will prove,   for $\mu$-a.e. $\xi\in \C$, that  the equation \eqref{E1} has a  unique  weak solution  with distribution $\P^\xi$ satisfying
\beq\label{*X} \P^\xi\bigg(\bigg\{\gg\in C([0,\infty);\C):\ \int_0^T |Z(\gg_s)|^2\d s<\infty, T>0\bigg\}\bigg)=1.\end{equation}  When $\tau=0$ and $P_t^0$ satisfies the Harnack inequality \eqref{G*2} below, the assertion holds for all initial point $x\in\H$ (in this case $\C=\H$).

To formulate the associated Markov semigroup, we introduce the process
\beq\label{RX} R^\xi(t):= \exp\bigg[  \int_0^t\< Z(X_s^\xi),\d W(s)\>- \ff 1 2 \int_0^t |Z|^2(X_s^\xi) \d s\bigg],\ \ t\ge 0,\end{equation} where $X_t^\xi$ is the segment solution to \eqref{E0} with inital value $\xi\in\C$.
By $\mu(\e^{\vv|Z|^2})<\infty$, this process is well defined for $\mu$-a.e. $\xi$.  We will use $S_t^Z$ to denote the semigroup of segment solutions to \eqref{E1}.
But when $\tau=0$, we use  $P_t^Z$ to replace $S_t^Z$  for  the notation  consistency  with Theorem \ref{T1.1}.

\beg{thm}\label{T2.1} Assume {\bf (A)} and   $ \mu(\e^{\vv|Z|^2})<\infty$ for some constant $\vv>0$.
\beg{enumerate}
\item[$(1)$] For $\mu$-a.s. $\xi\in \C$, $R^\xi(t)$ is a martingale and the equation $\eqref{E1}$ has a  weak solution starting at $\xi$ satisfying $\eqref{*X}$. Moreover, the associated  Markov semigroup $S_t^Z$ is given by
\beq\label{SM} (S_t^Zf)(\xi)= \E\big[f(X_t^\xi)R^\xi(t)\big],\ \ t\ge 0,\ \ f\in L^\infty(\mu). \end{equation}

\item[$(2)$] For any $\xi\in\C$, $\eqref{E1}$ has at most one weak solution with distribution $\P^\xi$ satisfying $\eqref{*X}$.
\item[$(3)$] Let $\tau=0$.  If there exist $p>1$ and $\Phi_p\in C((0,\infty)\times\H^2)$ such that
\beq\label{G*1} \int_0^t \ff{\d s}{\{\mu_0(\exp[-\Phi_p(s,x,\cdot)])\}^{\ff 1 p}}<\infty,\ \ t>0, x\in\H,\end{equation} and $P_t^0$ satisfies the Harnack inequality
\beq\label{G*2} (P_t^0f)^p(x)\le \e^{\Phi_p(t,x,y)}P_t^0f^p(y),\ \ t>0, f\in \B^+(\H), x,y\in\H,\end{equation}
then  for any $x\in \H$, $R^x(t)$ is a martingale,
 the equation $\eqref{E1}$ has a unique weak solution satisfying $\eqref{*X}$ starting at $x$, and  the Markov semigroup is given by
\beq\label{SM'}(P_t^Zf)(x)= \E\big[f(X^x(t))R^x(t)\big],\ \ t\ge 0,\ x\in\H, \ f\in \B_b(\H).\end{equation}\end{enumerate}\end{thm}

\paragraph{Remark 2.1.} Although the reference measure $\mu$ is less explicit, the condition $\mu(\e^{\vv|Z|^2})<\infty$ can be verified by using the marginal distribution $\mu_0$, which is explicitly given in applications, for instance, in \eqref{SDE} $\mu_0(\d x)=\e^{-V(x)}\d x$ for $\H=\R^d$. Let, for instance, $$Z(\xi)= \int_{-\tau}^0 h(\xi(\theta))\d\theta,\ \ \xi\in \C,$$ where $h$ is a measurable function on $\H$. Then
\beg{align*} \mu(\e^{\vv |Z|^2}) &= \E \e^{\vv |Z|^2(X_\tau^\mu)} = \E \e^{\vv\big|\int_{-\tau}^{0} h(X^{\mu_0}(\tau+\theta))\d\theta\big|^2}\\
&\le \E\e^{\vv \tau\int_0^\tau |h(X^{\mu_0}(s))|^2\d s} \le \ff 1 \tau\int_0^\tau \E \e^{\vv \tau^2 |h(X^{\mu_0}(s))|^2}\d s= \mu_0(\e^{\vv \tau^2 |h|^2}).\end{align*} Therefore, $\mu(\e^{\vv |Z|^2})<\infty$ follows from $\mu_0(\e^{\vv\tau^2h^2})<\infty.$

\

To prove Theorem \ref{T2.1}, we first present the following lemma.

\beg{lem} \label{L3} Assume {\bf (A)} and  $\mu(\e^{\vv|Z|^2})<\infty$ for some constant $\vv>0$.\beg{enumerate} \item[$(1)$] The process $$R^{\mu}(t):= \exp\bigg[  \int_0^t\<Z(X_s^{\mu}),\d W(s)\>- \ff 1 2 \int_0^t |Z|^2(X_s^{\mu}) \d s\bigg],\ \ t\ge 0$$ is a martingale. Consequently, for $\mu$-a.e. $\xi\in \C$, $R^\xi(t)$ is a martingale.

\item[$(2)$] Let   $\tau=0$. If  $\eqref{G*2}$ holds for some $p>1$ and $\Phi_p$ satisfying $\eqref{G*1}$,
then $R^x(t)$ is a martingale for any $x\in \H$. \end{enumerate} \end{lem}

\beg{proof}  (1) By the stationarity of $X_t^\mu$ and $\mu(\e^{\vv|Z|^2})<\infty$, we have
$$\E\int_0^T   |Z(X_s^\mu)|^2\d s = T\mu(|Z|^2) <\infty,\ \ T>0.$$ So, $\P$-a.s.
$$\tau_n  := \inf\bigg\{t\ge 0: \int_0^t |Z(X^\mu_s)|^2\d s \ge n\bigg\}\uparrow \infty\ \ {\rm as}\ n\uparrow\infty.$$  By Girsanov's  theorem (see e.g. \cite{L}), $(R^\mu(t\land\tau_n))_{t\ge 0}$ is a martingale for every $n\ge 1$, so that
 Fatou's lemma gives
\beg{align*} \E  (R^\mu(t)|\F_s)&= \E\Big(\liminf_{n\to\infty} R^\mu(t\land \tau_n)\Big|\F_s\Big)\le \liminf_{n\to\infty}\E\Big(  R^\mu(t\land \tau_n)\Big|\F_s\Big) \\
&= \liminf_{n\to\infty}R^\mu(s\land \tau_n)= R^\mu(s),\ \ t\ge s\ge 0.\end{align*} Thus, $(R^\mu(t))_{t\ge 0}$ is a supmartingale. Since $ \E\int_0^T   |Z(X_s^\mu)|^2\d s <\infty$ implies
$ \E\int_0^T  |Z(X_s^\xi)|^2\d s <\infty$ for   $\mu$-a.e. $\xi$, the above argument also implies that $(R^\xi(t))_{t\ge 0}$ is a supmartingale for $\mu$-a.e. $\xi$. Noting that $\E R^\mu(t)= \int_\C \E R^\xi(t) \d\mu$, we conclude that   $(R^\mu(t))_{t\ge 0} $ and   $(R^\xi(t))_{t\ge 0}$ for $\mu$-a.e. $\xi$ are martingales provided   $\E R^\mu(t)=1$ for all $t\ge 0$.

By the stationarity  of $X_t^\mu$ and  Jensen's inequality,  we have
$$\E \e^{\ff 1 2 \int_0^{2\vv} |Z|^2(X_s^{\mu})\d s} \le \ff 1 {2\vv} \int_0^{2\vv} \E \e^{\vv| Z|^2(X_s^{\mu})}\d s= \mu(\e^{\vv|Z|^2})<\infty.$$ Then Girsanov's theorem ensures  that
$(R^{\mu}(t))_{t\in [0, 2\vv]}$ is a martingale. In particular,
\beq\label{T01}\E R^\mu(t)=\E R^\xi(t)=1,\ \ t\in [0, 2\vv], \mu\text{-a.e.}\ \xi.\end{equation}
Assuming that   $\E R^{\mu}(t) =1 $ for $t\in [0, 2k\vv]$ and some $k\ge 1$,   it remains to prove   $\E R^{\mu}(t)=1$ for $t\in [2k\vv, 2(k+1)\vv].$
Let $t_1= 2\kk\vv$ and
$$\tt W(t)= W(t+t_1)-W(t_1),\ \ t\ge 0.$$ Then $\tt W(t)$ is a cylindrical Brownian motion on the same probability space with respect to filtration $\F_{t_1+t}.$ By {\bf (A)}, \eqref{E0} with $\tt W(t)$ replacing $W(t)$ has existence and uniqueness as well.   Let $\tt X^\xi_t$ be the segment process of the solution with  $\tt X_0^\xi= \xi$ defined as in \eqref{D*E}.    By the Markov property,
\beg{equation*}\beg{split}  \Xi(t) &:=  \E \big(\e^{\int_{t_1}^{t} \<Z(X_s^{\mu}), \d W(s)\>-\ff 1 2 \int_{t_1}^{t} |Z|^2(X_s^{\mu})\d s}|\F_{t_1}\big)\\
&=  \Big\{\E \big(\e^{\int_{0}^{t -t_1} \<Z(\tt X_s^{\xi}), \d \tt W(s)\>-\ff 1 2 \int_{0}^{t-t_1} |Z|^2(\tt X_s^\xi)\d s}\big)\Big\}\Big|_{\xi= X_{t_1}^{\mu}}\\
&= \big\{\E R^\xi(t -t_1)\big\}\big|_{\xi= X_{t_1}^{\mu}},\ \ t\ge t_1.\end{split}\end{equation*}
Since the law of $X_{t_1}^{\mu}$ is $\mu$, this and \eqref{T01} imply $\Xi(t)=1$ a.s. for all $t\in [t_1, t_1+2\vv]$, so that
$$\E(R^{\mu}(t)|\F_{t_1})= R^{\mu}(t_1)\Xi(t)=R^{\mu}(t_1).$$ So, by the assumption $\E R^\mu(t_1)=1$, we obtain  $\E R^{\mu}(t)  =1$ for $t\in [2k\vv, 2(k+1)\vv].$ 

(2) By  Girsanov's theorem and the Markov property, it suffices to find out a constant  $t>0$ such that
\beq\label{G*3} \E \e^{\ff 1 2 \int_0^{t}|Z(X^x(s))|^2\d s} <\infty,\ \ x\in\H.\end{equation} By \eqref{G*2}, we have
$$\mu_0\big(\e^{-\Phi_p(s,x,\cdot)}\big) \big(P_s^0 \e^{\ff\vv p|Z|^2}\big)^p (x) \le \mu_0\big(\e^{\vv |Z|^2}\big)<\infty.$$ Combining this with \eqref{G*1} and Jensen's inequality, for any $\ll>0$ and $t_\ll:= \ff{\vv}{p\ll}$ we obtain
\beq\label{G*4} \beg{split} & \E\big[\e^{\ll\int_0^{t_\ll} |Z(X^x(s) )|^2 \d s }\big]\le \ff 1 {t_\ll} \int_0^{t_\ll} P_s^0 \e^{\ff \vv p|Z|^2} (x)\d s \\
&\le \ff 1 {t_\ll} \int_0^{t_\ll} \big\{\mu_0(\e^{\vv |Z|^2})\big\}^{\ff 1 p} \big\{\mu_0(\e^{-\Phi_p(s,x,\cdot)})\big\}^{-\ff 1 p} \d s<\infty,\ \ x\in\H.\end{split}\end{equation}  In particular,   \eqref{G*3} holds for some constant $t>0$.
\end{proof}

\beg{proof}[Proof of Theorem \ref{T2.1}]
(1) By Lemma \ref{L3} and Girsanov's theorem, for $\mu$-a.s. $\xi$,
$$W^\xi(t):= W(t)-\int_0^tZ(X_s^\xi)\d s,\ \ t\in [0,T]$$ is a  cylindrical   Brownian motion on $\tt\H$ under the  probability measure $\Q^\xi$ defined on $\F_\infty$   by
$$\Q^\xi(A)= \E[1_A R^\xi(T)],\ \ T>0, A\in \F_T,$$ and $(X^\xi(t))_{t\ge 0}$ is a weak solution of \eqref{E1} with respect to the cylindrical Brownian motion $W^\xi(t)$. Note that $\Q^\xi$ is well defined according to the martingale property of $R^\xi(t)$ and the Kolmogorov consistency theorem.  Therefore, the associated  Markov   semigroup of the weak solution  is given by \eqref{SM}.

(2) For $\xi\in\C$  and  each $i=1,2$, let  $(X^{(i)}(t))_{t\ge 0}$ be a weak solution to \eqref{E1} starting at $\xi$  with respect to the cylindrical Brownian motion $W^{(i)}(t)$ under  a  complete filtration probability space  $(\OO^{(i)}, \{\F_t^{(i)}\}_{t\ge 0}, \P^{(i)}),$   such that the distribution  $\P_i^\xi$  satisfies   \eqref{*X}.  We  intend to prove $\P_1^\xi=\P_2^\xi$. By \eqref{*X}, we have
$$\tau^{(i)}_n := \inf\bigg\{t\ge 0: \int_0^t |Z(X^{(i)}_s)|^2 \d s  \ge n\bigg\}\uparrow\infty\ \text{as}\ n\uparrow\infty,\ \ i=1,2.$$   For every $i=1,2$ and $n\ge 1$,
$$R_n^{(i)}(t):= \exp\bigg[-\int_0^{\tau_n^{(i)}\land t} \big\<Z(X_s^{(i)}), \d W^{(i)}(s)\big\>- \ff 1 2 \int_0^{\tau_n^{(i)}\land t} |Z|^2(X_s^{(i)})\d s\bigg],\ \ t\ge 0$$ is a $\P^{(i)}$-martingale.
Define the probability  measure $\Q^{(i)}_n$ on $\F^{(i)}_\infty$   by letting
$$\Q_n^{(i)} (A)= \E_{\P^{(i)}} [1_A R_n^{(i)}(T)],\ \ T>0, A\in \F_T^{(i)}.$$  By  Girsanov's theorem,
$$\hat W^{(i)}(t):= W^{(i)}(t)+ \int_0^{t\land\tau_n^{(i)}} Z (X_s^\xi)\d s,\ \ t\ge 0$$ is a $\Q_n^{(i)}$-cylindrical  Brwonian motion   on $\tt\H$.  Therefore, up to time $\tau_n^{(i)}$,
$X^{(i)}(t)$ solves  the SDE \eqref{E0} with  the $\Q_n^{(i)}$-cylindrical Brownian motion $\hat W^{(i)}(t)$ replacing $W(t)$. By the pathwise (also  weak) uniqueness of \eqref{E0} according to {\bf (A)},  $(X^{(i)}(t), \hat W^{(i)}(t))_{t\in [0, T\land\tau_n^{(i)}]}$  under $\Q^{(i)}_n$ coincides in law with  $(X^\xi(t), W(t))_{t\in [0, T\land \tau_n^\xi]}$ under $\P$, where $T>0$ and
$$\tau_n^\xi:=\inf\bigg\{t\ge 0:\int_0^t |Z(X_s^\xi)|^2\d s\ge n\bigg\}.$$
 Therefore, for any $F\in \B_b(C([0,T];\H)\times C([0,T];\H))$,
\beg{equation*}\beg{split} &\E_{\P^{(i)}}\Big[1_{\{\tau_n^{(i)}\ge T\}}F\big(X^{(i)}([0,T]),  W^{(i)}([0,T])  \big)
\Big]\\
=   &\E_{\Q^{(i)}_n}\bigg[1_{\{\int_0^T|Z(X_s^{(i)}|^2 \d s \le n\}}\e^{\int_0^T \<Z(X_s^{(i)}), \d \hat W^{(i)}(s)\>- \ff 1 2 \int_0^{T} |Z|^2(X_s^{(i)})\d s} \\
&\qquad\ \ \times F\Big(X^{(i)}([0,T]), \Big(\hat W^{(i)} -\int_0^\cdot Z(X_s^{(i)})\d s\Big)([0,T])\Big)\bigg]
 \\
 &= \E \bigg[1_{\{\int_0^T|Z(X_s^\xi|^2 \d s\le n\}}\e^{\int_0^T \< Z(X_s^\xi), \d W(s)\>- \ff 1 2 \int_0^{T} |Z|^2(X_s^\xi)\d s} \\
&\qquad\ \ \times F\Big(X^\xi([0,T]), \Big(W-\int_0^\cdot Z(X_s^\xi)\d s\Big)([0,T])\Big)\bigg],\ \ i=1,2.\end{split} \end{equation*}
Consequently,
$$\E_{\P^{(1)}}\Big[1_{\{\tau_n^{(1)}\ge T\}} F\big(X^{(1)}([0,T]),  W^{(1)}([0,T])  \big)\Big] = \E_{\P^{(2)}}\Big[1_{\{\tau_n^{(2)}\ge T\}} F\big(X^{(2)}([0,T]),  W^{(2)}([0,T])  \big)\Big] $$ holds for any $n\ge 1.$
 Letting $n\to\infty$ we obtain
 $$\E_{\P^{(1)}}\Big[ F\big(X^{(1)}([0,T]),  W^{(1)}([0,T])  \big)\Big] = \E_{\P^{(2)}}\Big[ F\big(X^{(2)}([0,T]),  W^{(2)}([0,T])  \big)\Big] $$ for any $T>0$ and $F\in \B_b(C([0,T];\H)\times C([0,T];\H))$. Therefore,  $\P_1^\xi=\P_2^\xi$.

 (3) Let $\tau=0$. By Lemma \ref{L3}(2) and the Girsanov theorem,
$(X^x(t))_{t\ge 0}$ is a weak solution to \eqref{E1} satisfying \eqref{*X} for the $\tt\H$-cylindrical Brownian motion $$W^x(t):= W(t)-\int_0^t Z(X^x(s))\d s,\ \ t\ge 0 $$    under the probability measure $\Q^x$, which is defined on $\F_\infty$ by
$$\Q^x(A):=\E[1_A R^x(T)],\ \ T>0, A\in \F_T.$$
Then the proof is finished by combining this with (2).
\end{proof}

\section{Uniqueness of invariant probability measure}

Since by Theorem \ref{T2.1}  $S_t^Z$ is  a Markov semigroup  on $L^\infty(\mu)$, it is meaningful to consider the class of   invariant probability measures absolutely continuous with respect to $\mu$:
$$\scr P_Z:=\big\{\rr\mu:\ \rr\ge 0, \mu(\rr)=1, \mu(\rr S_t^Z f)=\mu(\rr f)\ \text{for}\ t\ge 0, f\in \B_b(\C)\big\}.$$
Recall that when $\tau=0$ we use $P_t^Z$ to replace $S_t^Z$.

\beg{thm}\label{T0} Assume {\bf (A)} and  $\mu(\e^{\vv|Z|^2})<\infty$ for some $\vv>0$.  \beg{enumerate} \item[$(1)$]  If there exists $t>0$ such that $P_t^0$ has a strictly positive density $p_t^0(x,y)$ with respect to $\mu_0$, then $\nu\in\scr P_Z$ implies that $\rr:=\ff{\d\nu}{\d\mu}$ has a strictly positive version, and $\scr P_Z$ contains at most one element.
\item[$(2)$] In the situation of Theorem $\ref{T2.1}(3)$, the Markov semigroup $P_t^Z$ defined on $\B_b(\H)$  has at most one invariant probability measure. \end{enumerate}
 \end{thm}

\beg{proof}  (1) Let $\nu=\rr\mu \in \scr P_Z$. We first prove that $\rr$ has a strictly positive version;  i.e. $\mu$ is absolutely continuous with respect to $\nu$.  For measurable $A\subset \C$ with $\nu(A)=0$, we intend to prove $\mu(A)=0$.  Since $\nu$ is $S_t^Z$-invariant, we have
$$ \int_\C \E[R^\xi(t +\tau)1_A(X_{t+\tau}^\xi)]\nu(\d \xi)  =\nu(S_t^Z1_A)=\nu(A)=0.$$ Noting that  $R^\xi(t+\tau)>0$ for $\mu$-a.e. (hence, $\nu$-a.e.) $\xi$, this implies
$$\nu(S_{t+\tau}^0 1_A)= \int_\C \E [1_A(X_{t+\tau}^\xi)]\nu(\d\xi)=0.$$
Letting  $\rr_0(x)=\mu(\rr|\xi(0)=x)$ be the regular conditional expectation of $\rr$ with respect to $\mu$ given $\xi(0)$, from this and   \eqref{3*} we obtain
$$\mu_0((\{P_t^0\}^* \rr_0)S_\tau 1_A)= \mu_0(\rr_0 P_t^0(S_\tau 1_A)) = \nu(S_{t+\tau}^0 1_A) =0,$$ where due to  $p_t^0>0$ and $\mu_0(\rr_0)=1$,
$$(P_t^0)^*\rr_0:=\int_\H p_t^0(z,\cdot)\rr_0(z)\mu_0(\d z)>0.$$ So, $\mu_0(S_\tau 1_A)=0.$ Combining this with   \eqref{3*} and that $\mu$ is $S_\tau^0$-invariant, we obtain $\mu(A)=\mu(S_\tau^0 1_A)=\mu_0(S_\tau 1_A)=0.$

Next, according to \cite[Proof of Proposition 3.1(3)]{WY11}, the uniqueness follows if $S_{t+\tau}^Z$ has a strictly positive density with respect to $\nu$. Since $\mu$ is equivalent to $\nu$ as proved above, and $S_{t+\tau}^Z$ is equivalent to $S_{t+\tau}^0$ according to \eqref{SM}, it suffices to prove that $S_{t+\tau}^0$ has a strictly positive density with respect to $\mu$. Let $(S_\tau^0)^*$ be the adjoint operator of $S_\tau^0$ in $L^2(\mu)$, and let $\hat p_t^0(\xi,\eta)= p_t^0(\xi(0), \eta(0))$.
For any $f\in \B_b(\C)$,  \eqref{3*} yields
\beg{align*}& S_{t+\tau}^0 f(\xi)=P_t^0(S_\tau f)(\xi(0)) = \int_\H (S_\tau f)(y) p_t^0(\xi(0),y)\mu_0(\d y)\\
&= \int_\C (S_\tau^0 f)(\eta)  \hat p_t^0(\xi, \eta) \mu(\d\eta)= \int_\C f(\eta) (S_\tau^0)^* \hat p_t^0(\xi,\cdot)(\eta) \mu(\d\eta). \end{align*}
Since $p_t^0>0$ implies $\hat p_t^0>0$,  this implies that $S_{t+\tau}^0$ has a strictly positive density $(\xi,\eta)\mapsto (S_\tau^0)^*\hat p_t^0(\xi,\cdot)(\eta).$

(2)  By \cite[Proposition 3.1]{WY11}, $P_t^Z$ has at most one invariant probability measure if there exist  $t>0, q>1$ and a measurable function $\Psi: \H^2\to (0,\infty)$ such that
\beq\label{HS} P^Z_tf(x)\le (P^Z_tf^q(y))^{\ff 1 q}\Psi(x,y),\ \ f\in \B^+(\H), x,y\in \H.\end{equation}
By \eqref{G*4},   for any $r>1$ there exists a constant $t(r)>0$ such that
$$\GG_r(x,t):= \E \big[(R^x(t))^r+ (R^x(t))^{-r}\big]<\infty,\ \ t\in [0, t(r)], x\in \H.$$ Then \eqref{SM'} and \eqref{G*2} yield
\beg{align*} &P_t^Zf(x)= \E\big[f(X^x(t)R^x(t)\big] \le\ss{\GG_2(x,t) P_t^0 f^2 (x)}\le \ss{\GG_2(x,t)} \big(P_t^0 f^{2p}(y))^{\ff 1 {2p}}\e^{\ff 1{2p}\Phi_p(t,x,y)}\\
&\le \ss{\GG_2(x,t)} \big(\E [f^{4p}(X^y(t))R^y(t)]\big)^{\ff 1 {4p}}\big(\E[(R^y(t))^{-1}]\big)^{\ff 1 {4p}} \e^{\ff 1{2p}\Phi_p(t,x,y)}\\
&\le \ss{\GG_2(x,t)} \big(\GG_1(y,t)\big)^{\ff 1 {4p}} \e^{\ff 1{2p}\Phi_p(t,x,y)} \big(P_t^Z   f^{4p}(y))^{\ff 1 {4p}}.\end{align*} Therefore, \eqref{HS} holds for $q=4p$ and some function $\Psi$, and  the proof is thus finished.
\end{proof}

\section{Entropy estimate   using hyperboundedness}

In this section,  we assume that $P_t^0$ is hyperbounded, i.e. there exist $p_0>1$ and $t_0>0$ such that
\beq\label{HP} \|P_{t_0}^0\|_{L^2(\mu_0)\to L^{2p_0}(\mu_0)}:= \sup_{\mu_0(f^2)\le 1} \mu_0\big(|P_{t_0}^0f|^{2p_0}\big)^{\ff 1 {2p_0}}<\infty.\end{equation}According to Gross \cite{Gross},     when $P_t^0$ is symmetric in $L^2(\mu_0)$, for instance, $Z_0$ is given by $\eqref{Z0}$ and $\mu_0(\d x)= \e^{-V(x)}\d x$,    \eqref{HP} is equivalent to the defective log-Sobolev inequality \eqref{LS} for some constants $\kk>0, \bb\ge 0$. However, in the non-symmetric case, the latter is strictly  stronger than \eqref{HP}, see Examples 4.1 and 4.2 below for hypercontractive Markov semigroups without the log-Sobolev inequality. So, the following result  is new even for $\tau=0$ and $\H=\R^d$.

\beg{thm}\label{T2.2} Assume {\bf (A)} and $\eqref{HP}$ for some $t_0>0$ and $p_0>1$. If $  \mu(\e^{\ll |Z|^2})<\infty$ for some constant $\ll> \ff{(3p_0-1)(t_0+\tau)}{2(p_0-1)}$, then there exists $\nu:=\rr\mu\in \scr P_Z$
 such that
\beq\label{ENT} \mu(\rr\log\rr)\le \ff{(t_0+\tau) (3p_0-1)\log  \mu(\e^{\ll |Z|^2}) + 4\ll p_0 \log\|P_{t_0}^0\|_{L^2(\mu_0)\to L^{2p_0}(\mu_0)} }{2\ll(p_0-1)- (3p_0-1)(t_0+\tau)}. \end{equation}
\end{thm}
\beg{proof} Let $c=t_0+\tau$ and
 \beq\label{SEQ}\nu_n:= \ff 1 {cn} \int_0^{cn} \mu S_t^Z \d t,\ \ n\ge 1,\end{equation} where the probability measure $\mu S_t^Z$ is defined by
 $(\mu S_t^Z)(A):= \mu(S_t^Z 1_A), A\in \B(\C).$  It suffices to find a subsequence $n_k\to \infty$ such that $\nu_{n_k}\to \nu$ weakly
 for some probability measure $\nu:=\rr\mu$ with density $\rr$ satisfying \eqref{ENT}. We complete the proof by the following three steps.

 (a) Let $c_0=\|P_{t_0}^0\|_{L^2(\mu_0)\to L^{2p_0}(\mu_0)}<\infty$ and $p=1+\ff{p_0-1}{2p_0}>1$.  We first prove
 \beq\label{T1} \E\big[\e^{\int_0^{cn} f(X_s^\mu)\d s}\big]\le c_0^{\ff n{p}} \Big\{\mu\big(\e^{\ff{c(3p_0-1)}{p_0-1}f}\big)\Big\}^{\ff{n(p_0-1)}{2p_0p}},\ \ f\in \B_b(\C). \end{equation}
 Since   $\mu$ is an invariant probability measure of the segment process $X_t$, \eqref{3*}  implies
$$\mu_0(S_{t+\tau} f)= \int_{\C} \E[ f( X_{t+\tau}^{\xi(0)})]\mu(\d \xi)= \int_{\C} \E[ f( X_{t+\tau}^{\xi})]\mu(\d \xi)= \E [f(X_{t+\tau}^{\mu}) ]  =\mu(f)$$ for $t\ge 0$ and $ f\in \B_b(\C).$ Combining this with \eqref{3*} and using Jensen's inequality,
we obtain
\beq\label{3*1} \beg{split} &\mu_0(|S_c f|^{2p_0}) = \mu(|P_{t_0}^0(S_\tau f)|^{2p_0}) \le c_0^{2p_0} \mu_0((S_\tau f)^2)^{p_0}\\
&\le c_0^{2p_0}\big\{\mu_0(S_\tau f^2)\big\}^{p_0} = c_0^{2p_0} \big\{\mu(f^2)\big\}^{p_0},\ \ \ f\in L^2(\mu).\end{split}\end{equation}
For any $F\in \B_b(\C)$,   we consider  the Feymann-Kac semigroup
$$(\bar S_t^Ff)(\xi):= \E\Big[f(X_t^\xi)\e^{\int_0^t F(X_s^\xi)\d s}\Big],\ \ f\in \B_b(\C). $$
By the H\"older/Jensen inequalities, \eqref{3*} and \eqref{3*1}, we obtain
\beg{align*} &\mu(|\bar S_{c}^Ff|^{2p}) =\int_\C \Big(\E\big[f(X_{c}^\xi)\e^{\int_0^{c} F(X_s^\xi)\d s}\big]\Big)^{2p}\mu(\d \xi)\\
&\le \int_\C \Big\{\big(\E [f^{p}(X_{c}^\xi)]\big)^2\big(\E\e^{\ff{p}{p-1}\int_0^{c} F(X_s^\xi)\d s}\big)^{2(p-1)}\Big\} \mu(\d\xi)\\
&\le \int_\C  \big(S_{c}f^{p}\big)^2(\xi(0)) \bigg(\ff 1 {c} \int_0^{c} \E \e^{\ff{cp}{p-1}F(X_s^\xi)} \d s\bigg)^{2 (p-1)} \mu(\d\xi)\\
&\le \bigg(\int_{\H}   \big(S_{c}f^{p}\big)^{2p_0} \d\mu_0\bigg)^{\ff 1 {p_0}}  \bigg\{\int_\C \bigg(\ff 1 {c} \int_0^{c}   \E \big[\e^{\ff{cp}{p-1}F(X_s^\xi)}\big] \d s \bigg)^{\ff{2p_0(p-1)}{p_0-1}}\mu(\d\xi)\bigg\}^{\ff {p_0-1}{p_0}}\\
&\le c_0^{2} \mu(f^{2p}) \Big\{\mu\big(\e^{\ff{c(3p_0-1)}{p_0-1}F}\big)\Big\}^{\ff{p_0-1}{p_0}}<\infty.\end{align*}
So, $\bar S_{c}^F$ is bounded in $L^{2p}(\mu)$ with
$$\|\bar S_{c}^F\|_{L^{2p}(\mu)} \le c_0^{\ff 1 {p}} \Big\{\mu\big(\e^{\ff{c(3p_0-1)}{p_0-1}F}\big)\Big\}^{\ff{p_0-1}{2p_0p}}.$$   By the semigroup property, this leads to
$$\E\e^{\int_0^{cn} F(X_s^{\mu})\d s}= \mu(\bar S_{cn}^F1) \le \|\bar S_{cn}^F\|_{L^{2p}(\mu)} \le \|\bar S_c^F\|_{L^{2p}(\mu)}^n\le c_0^{\ff n{p}} \Big\{\mu\big(\e^{\ff{c(3p_0-1)}{p_0-1}F}\big)\Big\}^{\ff{n(p_0-1)}{2p_0p}}.$$ Therefore, \eqref{T1} holds.

(b)   For any $\vv\in (0,1),$ we intend to prove
\beq\label{RM} \E\big[R^{\mu}(cn)\log R^{\mu}(cn)\big]\le \ff{n\vv}{p(1-\vv)}\log\Big[c_0\big\{\mu\big(\e^{\ff{c(3p_0-1)}{2\vv(p_0-1)}|Z|^2}\big)\big\}^{\ff{p_0-1}{2p_0}}\Big].\end{equation}
We will   apply the following Young inequality (see \cite[Lemma 2.4]{ATW09}): for a probability measure $\LL$ on a measurable space $M$,
\beq\label{Young} \LL(fg)\le \LL(f\log f) + \log \LL(\e^g),\ \ f,g\in \B^+(M), \LL(f)=1.\end{equation} By the Young  inequality and \eqref{T1},
\beg{align*} & \E\big[R^{\mu}(cn)\log R^{\mu}(cn)\big] = \ff 1 2 \E \bigg[R^{\mu}(cn)\int_0^{cn} |Z(X_s^{\mu})|^2\d s\bigg]\\
&\le\vv\E\big[R^{\mu}(cn)\log R^{\mu}(cn)\big] + \vv\log\E\big[\e^{\ff 1 {2\vv} \int_0^{cn} |Z(X_s^{\mu})|^2\d s}\big]\\
&\le  \vv\E\big[R^{\mu}(cn)\log R^{\mu}(cn)\big] + \ff{n\vv}{p}\log \Big[c_0\big\{\mu\big(\e^{\ff{c(3p_0-1)}{2\vv(p_0-1)}|Z|^2}\big)\big\}^{\ff{p_0-1}{2p_0}}\Big].\end{align*}
When $Z$ is bounded we have $\E\big[R^{\mu}(cn)\log R^{\mu}(cn)\big]\le \ff {cn}2 \|Z\|_\infty^2<\infty$, so that this implies \eqref{RM}.
In general, let $Z_m= Z1_{\{|Z|\le m\}}, m\ge 1.$ Define $R^\mu_m(t)$ as $R^\mu(t)$ using $Z_m$ replacing $Z$. Then the assertion for bounded $Z$ implies
$$\E\big[R_m^{\mu}(cn)\log R_m^{\mu}(cn)\big]\le \ff{n\vv}{p(1-\vv)}\log\Big[c_0\big\{\mu\big(\e^{\ff{c(3p_0-1)}{2\vv(p_0-1)}|Z|^2}\big)\big\}^{\ff{p_0-1}{2p_0}}\Big],\ \ m\ge 1.$$
Due to Fatou's lemma,  we prove \eqref{RM}  by letting $m\to\infty$.

(c) By \eqref{SM}, \eqref{3*1}, \eqref{RM} and \eqref{Young}, for any $f\ge 0$ we have
\beq\label{FF0} \beg{split}\nu_n(f) &= \ff 1 {cn} \int_0^{cn} \E[R^\mu(cn) f(X_s^\mu)]\d s \\
&\le \ff 1 {cn} \E\big[R^{\mu}(cn)\log R^{\mu}(cn)\big] + \ff 1 {cn} \log\E\big[\e^{\int_0^{cn} f(X_s^\mu)\d s}\big]\\
 &\le    \ff{\vv}{cp(1-\vv)}\log\Big[c_0\big\{\mu\big(\e^{\ff{c(3p_0-1)}{2\vv(p_0-1)}|Z|^2}\big)\big\}^{\ff{p_0-1}{2p_0}}\Big]  + \ff 1 {cp} \log \Big[c_0 \Big\{\mu\big(\e^{\ff{c(3p_0-1)}{p_0-1}f}\big)\Big\}^{\ff{p_0-1}{2p_0}}\Big].\end{split}\end{equation}
 This implies that $\nu_n$ is absolutely continuous with respect to $\mu$, and  \eqref{ENT} holds for $\rr_n:=\ff{\d\nu_n}{\d\mu}$ replacing $\rr$. Indeed,
taking $f= R1_A$ in \eqref{FF0} for $\mu(A)=0$ and $R>0$, we obtain
$$\nu_n(A)\le  \ff 1 R\bigg( \ff{\vv}{cp(1-\vv)}\log \Big[c_0\big\{\mu\big(\e^{\ff{c(3p_0-1)}{2\vv(p_0-1)}|Z|^2}\big)\big\}^{\ff{p_0-1}{2p_0}}\Big]+ \ff 1 {cp} \log c_0\bigg),\ \ R>0.$$
 Letting $R\to\infty$ we prove $\nu_n(A)=0$ for $\mu(A)=0$, so that $\nu_n$ has a density $ \rr_n$ with respect to $\mu$. Next,
applying \eqref{FF0}  to $f= \ff{p_0-1}{c(3p_0-1)}\log(\rr_n\land m +m^{-1})$ and letting $m\to\infty$, we obtain
 $$\ff{p_0-1}{c(3p_0-1)}\mu(\rr_n\log\rr_n) \le \ff{\vv}{cp(1-\vv)}\log \Big[c_0\big\{\mu\big(\e^{\ff{c(3p_0-1)}{2\vv(p_0-1)}|Z|^2}\big)\big\}^{\ff{p_0-1}{2p_0}}\Big]+ \ff 1 {cp}\log c_0.$$
 Taking $\vv= \ff{c(3p_0-1)}{2\ll (p_0-1)}\in (0,1)$ such that $\ff{c(3p_0-1)}{2\vv(p_0-1)}=1,$ we arrive at
 \beg{align*} \mu(\rr_n\log\rr_n) &\le \ff{\vv(3p_0-1)}{p(1-\vv)(p_0-1)} \log\big[c_0\{\mu(\e^{\ll|Z|^2})\}^{\ff{p_0-1}{2p_0}}\big] +\ff{3p_0-1}{p(p_0-1)}\log c_0\\
 &= \ff{2cp_0(3p_0-1)}{[2\ll(p_0-1)-c(3p_0-1)](p_0-1)} \log\big[c_0\{\mu(\e^{\ll|Z|^2})\}^{\ff{p_0-1}{2p_0}}\big] +\ff{2p_0}{p_0-1} \log c_0\\
 &= \ff{c (3p_0-1)\log  \mu(\e^{\ll |Z|^2}) + 4\ll p_0 \log c_0 }{2\ll(p_0-1)- c(3p_0-1)},\ \ n\ge 1.\end{align*}
  Therefore, $\{\rr_n\}_{n\ge 1}$ is uniformly integrable in $L^1(\mu)$, so that for some subsequence $n_k\to\infty$ we have $\rr_{n_k}\to \rr$ weakly in $L^1(\mu)$. Then $\nu_{n_k}\to \nu:=\rr\mu$  strongly  and $\rr$ satisfies \eqref{ENT}.
\end{proof}

We first consider a simple example to show that the integrability condition in Theorem \ref{T2.2} is asymptotically sharp for small $t_0$.

\paragraph{Example 4.1} Let $\H=\R^d$, $\tau=0, \si=\ss 2 I$ and $Z_0=-x$. Then $\mu_0(\d x)= (2\pi)^{-\d/2} \e^{-\ff 1 2|x|^2}\d x$ is the standard Gaussian measure. It is well known by Nelson \cite{Nel} (see also Gross \cite{Gross}), we have
$$\|P_{t_0}^0\|_{L^2(\mu_0)\to L^{2p_0}(\mu_0)}=1,\ \ t_0>0, p_0=\ff 1 2(1+\e^{2t_0}).$$
Since $\ff{3p_0-1)t_0}{2(p_0-1)}= \ff{t_0(3\e^{2t_0}-1)}{2(\e^{2t_0}-1)}\to \ff 1 2$ as $t_0\to 0$,  for any $\ll>\ff 1 2$  there exists $t_0>0$ such that $\ll> \ff{3p_0-1)t_0}{2(p_0-1)}.$ By Theorems \ref{T0} and \ref{T2.2}, if
$\int_{\R^d}\e^{\ll |Z(x)|^2- \ff 1 2|x|^2}\d x<\infty$ then there exists a unique $\nu=\rr\mu_0\in\scr P_Z$ and
$$\mu_0(\rr\log \rr)\le \inf\Big\{  \ff{t_0 (3\e^{2t_0}-1)\log \mu_0(\e^{\ll|Z|^2})}{2\ll(\e^{2t_0}-1)- t_0 (3\e^{2t_0}-1)}:  \ff{t_0(3\e^{2t_0}-1)}{2(\e^{2t_0}-1)}<\ll\Big\}<\infty. $$
On the other hand, for any $\ll<\ff 1 2$, there exists $Z$ with $\mu_0(\e^{\ll|Z|^2})<\infty$ but $P_t^Z$ does not have any invariant probability measure. For instance, when $Z(x)=x$ we have $\mu_0(\e^{\ll |Z|^2})<\infty$ for any $\ll<\ff 1 2$ but $P_t^Z=\e^{t\DD}$ does not have invariant probability measure.

\

Below we   consider   three more examples.  The first two are  degenerate SDEs and semilinear SPDEs for which the defective log-Sobolev inequality does not hold, and the last belongs to  monotone  SPDEs where the defective log-Sobolev inequality is unknown. See \cite{Liu, W07, W15} for more examples of hyperbounded  Markov semigroups without the defective log-Sobolev inequality.

\paragraph{Example 4.2 (Infinite-dimensional stochastic Hamiltonian system).} Let   $\H_1$ be a separable Hilbert space. Consider the following SPDE for $(X(t),Y(t))$ on $\H:=\H_1\times\H_1$:
$$ \beg{cases}\d X(t)=\{ Y(t)-LX(t)\}\d t,\\
\d Y(t)= \{Z(X_t,Y_t)-L Y(t)\}\d t + \d W(t),\end{cases}$$
where $W(t)$ is the standard cylindrical Brownian motion on $\H_1$, $Z: \C\mapsto \H_1$  is measurable, $(L,\D(L))$ is a positive definite self-adjoint operator on $\H_1$ with discrete eigenvalues $0<\ll_1\le\ll_2\le\cdots$ satisfying $\sum_{i=1}^\infty \ll_i^{-\dd}<\infty$ for some constant $\dd\in (0,1).$ Then the reference SDE with $Z=0$ has a unique invariant probability measure
$\mu_0= N(0, (2L)^{-1})\times N(0, (2L)^{-1})$, where  $N(0, (2L)^{-1})$ is the centered Gaussian measure on $\H_1$ with covariance $(2L)^{-1}.$ By \cite[Theorem 4.1]{W15b} with $Z=0, A=0, B=I, L_1=L_2=L$ such that $\dd=0$, the associated Markov semigroup $P_t^0$ is
hypercontractive. So,    Theorem \ref{T2.2} applies. Moreover, the Harnack inequality in \cite[Lemma 4.2]{W15b} implies that $P_t^0$ has a strictly positive density with respect to $\mu_0$. Then  Theorem \ref{T0} implies the uniqueness of invariant probability measure of $S_t^Z$.

\paragraph{Example 4.3 (Finite-dimensional stochastic Hamiltonian system).} Consider the following degenerate SDE for $(X(t),Y(t))$ on $\H= \R^{2d}$:
$$ \beg{cases}\d X(t)=  Y(t) \d t,\\
\d Y(t)= \{Z(X_t,Y_t)- X(t)- Y(t)\}\d t +   \d W(t),\end{cases}$$ Let $P_t^0$ be the Markov semigroup for the SDE with $Z=0$.  By  \cite[Theorem 4.4]{GW},    for any $p>1$ there exists a constant $c>0$ such that
$$(P_tf)^p(x)\le (P_tf^p)(y) \e^{\ff{c|x-y|^2}{1\land t^3}},\ \ t>0, x,y\in \R^{2d}, f\in \B^+(\R^{2d}).$$   Since $\Phi_p(s,x,y):= \ff{c|x-y|^2}{1\land s^3}$ satisfies
$$\mu_0(\e^{-\Phi_p(s,x,\cdot)})\ge \e^{-c} \mu_0(B(x, 1\land s^{\ff 3 2}))\ge \aa(x) (1\land s)^{\ff{3 d}2},\ \ s>0, x\in \R^{2d} $$ for some positive $\aa\in C(\R^{2d})$, \eqref{G*1} holds for $p>\ff{3 d} 2.$ Therefore, when $\mu(\e^{\vv|Z|^2})<\infty$ for some $\vv>0$,   Theorem \ref{T0}(2) implies that $P_t^Z$ is a Markov semigroup on $\B_b(\R^{2d})$ having at most one invariant probability measure.

Moreover,   by \cite[Example 5.1]{W15b}   with $W=0$,    $P_t^0$ has  unique invariant probability measure $\mu_0(\d x):= (2\pi)^{-d } \e^{-\ff 1 2 |x|^2}\d x$ on $\R^{2d} $, and
 $$\|P_{t_0}^0\|_{L^2(\mu_0)\to L^4(\mu_0)}=1$$ holds for some constant $t_0>0$. Therefore, by Theorem \ref{T2.2}, if $\mu(\e^{\ll|Z|^2})<\infty$ for some $\ll> 4 t_0$ then $P_t^Z$ has a (unique, as observed above) invariant probability measure $\nu$ with density $\rr:=\ff{\d\nu}{\d\mu_0}$ satisfying
 $$\mu_0(\rr\log \rr) \le \ff{4 t_0 \log\mu(\e^{\ll |Z|^2})}{\ll - 4 t_0}. $$

 \paragraph{Example 4.4 (Monotone SPDE).} Let $\si(0)\in \scr L_{HS}(\tt\H,\H)$, the class of Hilbert-Schmidt operators from $\tt\H$ to $\H$, and let $\aa\ge 2$ be a constant. Assume that
 $r\mapsto \,_{\V^*}\<Z_0(r v_1+v_2),v_3\>_\V$ is continuous in $r\in\R$ for any $v_1,v_2,v_3\in\V$, and there exist constants $C,\dd>0$ such that
 \beg{align*} & 2_{\V^*}\<Z_0(v_1)-Z_0(v_2), v_1-v_2\>_\V+ \|\si(v_1)-\si(v_2)\|_{L_{HS}}^2\le C \|v_1-v_2\|_\H^2-\dd \|v_1-v_2\|_\V^\aa,  \\
 &\|Z_0(v)\|_{\V^*}\le C+C\|v\|_\V^{\aa-1} \ \ {\rm for \ all\ } v_1,v_2,v\in\V.\end{align*}
According to \cite[Theorem 1.4]{Liu},   the equation \eqref{E0} has a unique variational solution for any initial value and   the Markov semigroup $P_t^0$ is hyperbounded with respect to the unique invariant probability measure $\mu_0$. Moreover, according to \cite{WY11}, the Harnack inequality in \cite[Theorem 1.2]{Liu} implies that $P_t^0$ has a positive density with respect to $\mu_0$. So, Theorem \ref{T0}(1)  and Theorem \ref{T2.2} apply. When  $\aa>2$,  \cite[Theorem 1.4]{Liu} ensures   \eqref{HP}  for any $t_0>0$ and $p_0>1>0$,   so that by these results   $\mu(\e^{\ll|Z|^2})<\infty$ for some $\ll>4\tau$ implies that $\scr P_Z$ contains a unique measure $\nu$, which has a strictly positive density $\rr$ with respect to $\mu$, and $\mu(\rr\log\rr)<\infty.$

\section{Sobolev estimates  using log-Sobolev inequality }

In this section, we aim to extend Theorem \ref{T1.1} to   degenerate SDEs with path-dependent drifts. When $\tau>0$, we will consider the Sobolev regularity of the marginal density  of the invariant probability measure.
  For a  probability measure $\nu$ on $\C$ and $\theta\in [-\tau,0]$, let $\nu_\theta$ be the $\theta$-marginal distribution of $\nu$, i.e.
$$\nu_\theta(A):= \nu(\{\xi\in \C: \xi(\theta)\in A\}),\ \ A\in\B(\H).$$ In particular, by the stationarity of $X^{\mu_0}(t)$, we have $\mu_\theta=\mu_0$ for $\theta\in [-\tau,0].$

We mainly consider the finite-dimensional case, but  make a simple infinite-dimensional extension in \S 5.4.
Let   $\H=\R^d$ and $\tt \H=\R^m$ for some $d,m\ge 1$, and let $V\in C^2(\R^d)$ such that $\mu_0(\d x):=\e^{-V(x)}\d x$ is a probability measure on $\R^d$.  Let  $\si\in C^2(\R^d\to \R^d\otimes\R^m) $ and let $Z_0$ be in \eqref{Z0}.  Then the operator
\beq\label{SL} \scr L_0:= \ff 1 2 \sum_{i,j=1}^d (\si\si^*)_{ij} \pp_i\pp_j + \sum_{i=1}^d \<Z_0, e_i\>\pp_i\end{equation} defined on $C_0^\infty(\R^d)$ is symmetric  in $L^2(\mu_0)$; namely,
$$-\mu_0(f\scr L_0g)= \EE_0(f,g):= \mu_0(\<\si^*\nn f,\si^*\nn g\>),\ \ \ f,g\in C_0^\infty(\R^d).$$  Let     $H_\si^{1,2}(\mu_0)$ be the closure of $C_0^\infty(\R^d)$ with respect to the Sobolev norm
$$\|f\|_{H_\si^{1,2}(\mu_0)}:=\{\mu_0(|f|^2+|\si^*\nn f|^2)\}^{\ff 1 2}.$$ Then   $(\EE_0, H_\si^{1,2}(\mu_0))$ is a symmetric Dirichlet form on $L^2(\mu_0)$ and the associated Markov process can be constructed as the solution to the SDE
\beq\label{EE1} \d X(t) = Z_0(X(t))\d t +\si(X(t))\d W(t), \end{equation} where $W(t)$ is the $m$-dimensional Brownian motion.

As in Section 4, we investigate    the following functional SDE
\beq\label{EE2} \d X(t) = \{\si(X(t))Z(X_t)+Z_0(X(t))\}\d t +\si(X(t))\d W(t) \end{equation} by using integrability conditions on the measurable map   $Z: \C\to \R^m.$ Besides the existence of invariant probability measure and the entropy estimate presented  in Theorem \ref{T2.2}, we aim to derive more and stronger density estimates including those given in Theorem \ref{T1.1}.  To this end, we make the following assumption {\bf (H)}, where the log-Sobolev inequality is essentially stronger than the hyperboundedness of $P_t^0$ used in Section 4.

Let
{\rm Lie}$\{U_1,\cdots, U_m\}$ be the Lie algebra induced by vector fields $\{U_1,\cdots, U_m\}$. More precisely, let $\scr A_0=\{U_i: 1\le i\le m\}$ and $$\scr A_k = \big\{[U,U']:= UU'-U'U:\ U,U'\in \cup_{0\le l\le k-1} \scr A_l\big\}, \ \ k\ge 1.$$ Then  {\rm Lie}$\{U_1,\cdots, U_m\}$ is the linear space spanned by $\cup_{k\ge 0} \scr A_k.$

\beg{enumerate} \item[{\bf(H)}]  \emph{Let $V\in C^2(\R^d)$ such that $\mu_0(\d x):=\e^{-V(x)}\d x$ is a probability measure on $\R^d$. There exists $k\ge 2$ such that $\si\in C^k(\R^d\to\R^d\otimes\R^m)$ and  vector fields  $$U_i:= \sum_{j=1}^d \si_{ji}\pp_j,\ \ 1\le i\le m $$ satisfy the H\"ormander condition up to the $k$-th order of Lie brackets. Let $Z_0$ be in \eqref{Z0}.
Moreover, $1\in H_\si^{1,2}(\mu_0)$ with $\EE_0(1,1)=0$,  and the log-Sobolev inequality \eqref{LS} holds   for some constants $\kk>0$ and $\bb\ge 0.$}\end{enumerate}

This assumption implies that the solution to \eqref{EE1} is non-explosive, the associated Markov semigroup $P_t^0$ has strictly positive symmetric density $p_t^0(x,y)$ with respect to the unique invariant probability measure $\mu_0$, and $P_t^0$ is hyperbounded. Since the Dirichlet form is irreducible, the hyperboundedness of $P_t^0$ is equivalent to $\|P_t^0\|_{L^2(\mu_0)\to L^4(\mu_0)}=1$ for large $t>0$ (i.e. the  hypercontractivity), see
\cite{M, W14}. Consequently, the Poincar\'e inequality
\beq\label{P} \mu_0(f^2)\le C\EE_0(f,f),\ \ f\in H_\si^{1,2}(\mu_0), \mu_0(f)=0\end{equation} holds for some constant $C>0$.

\subsection{Main result and examples}

Let $\tau\ge 0$ and let  $\mu$ on $\C:=C([-\tau,0];\R^d)$  be  the unique invariant probability measure of the segment solution to \eqref{EE1}.
 We will need the condition $\mu(\e^{\ll|Z|^2})<\infty$ for   $\ll$   big enough in terms of $\kk$ and   $\tau$.
 Let $$\ll_{\kk,\tau}=\inf\Big\{\ll>\kk:\ \Big(1+\ss{1+8\ll/\tau}\big)\big(1- \ss{\kk/\ll}\big)\ge 16\Big\}.$$ When $\tau=0$, we have $\ll_{\kk,0}=\kk.$
 Then for any $\ll>\ll_{\kk,\tau}$ we have
 $$q_\ll:= \ff{2\ss\ll}{\ss\ll-\ss\kk +\ss{(\ss\ll-\ss\kk)^2 -\ff{16\ss\ll(\ss\ll-\ss\kk)}{1+\ss{1+8\ll/\tau}}}}\ge \ff{\ss\ll}{\ss\ll-\ss\kk}>1.$$

\beg{thm}\label{T3.1} Assume {\bf (H)}.  If $\mu(\e^{\ll|Z|^2})<\infty$  for some $\ll>\ll_{\kk,\tau}$,   then $\scr P_Z$ contains a unique probability measure $\nu=\rr\mu$. Moreover:
\beg{enumerate} \item[$(1)$]  For any $p\in (1, \ff{q_\ll}{q_\ll-1})$, there exists a constant $k=k(p,\ll)\ge 1$ such that
 \beq\label{SI4}\mu(\rr^p)\le k\mu(\e^{\ll|Z|^2})^k.\end{equation}
 \item[$(2)$] The marginal density $\rr_\theta:=\ff{\d \nu_\theta}{\d\mu_0}$ does not depend on $\theta\in [-\tau,0]$, and has a continuous, strictly positive version such that  $\log\rr_0, \rr_0^{\ff p 2}\in H^{1,2}_\si(\mu_0)$ for $p\in (1,\ff{q_\ll}{q_\ll-1})$ with
the following estimates holding  for some constant $k=k(p,\ll)\ge 1$:
\beq\label{SI3} \beg{split} & \mu_0\big(|\si^* \nn \ss{\rr_0}|^2\big)\le \ff 1 {\ll-\kk} \big\{\log\mu(\e^{\ll|Z|^2}) + \bb\big\}<\infty;\\
&\mu_0\big(|\si^*\nn\log \rr_0|^2\big)\le 4 \mu(|Z|^2) <\infty;
\\& \mu_0\big(|\si^*\nn\rr_0^{\ff p 2}|^2+ \rr_0^p \big)\le k \big\{\mu (\e^{\ll|Z|^2})\big\}^k.\end{split}\end{equation}\end{enumerate}
\end{thm}

  Since  Theorem \ref{T1.1} does not imply $\rr^{\ff p 2}\in H^{1,2}_\si(\mu_0)$, the last estimate is new even in the non-degenerate case   without delay (i.e. $\tau=0$).
We present below two examples of degenerate diffusion processes satisfying the log-Sobolev inequality such that Theorem \ref{T3.1} applies.

 \paragraph{Example 5.1 (Gruschin type diffusions).} Let $d=2$ and $l\in\mathbb N.$ Let
 $$U_1(x)= \pp_1,\ \ \ \ U_2(x)= x_1^l \pp_2,\ \ x=(x_1,x_2)\in\R^2.$$ Then the H\"ormander condition is satisfied. Let $m\ge 2, c_1\in\R, c_2\ne 0$ and $c_3,c_4>0$ such that $\mu_0(\d x):=\e^{-V(x)}\d x$ is a probability measure for
 $$V(x):= c_1 + (c_3|x_1-c_1|^{l+1}+c_4 x_2^2)^m.$$ Obviously, $1\in H_\si^{1,2}(\mu_0)$ with $\EE_0(1,1)=0$. Moreover, by \cite[Proposition 4.1]{W09}, \eqref{LS} holds for some constant $\kk>0$ and $\bb=0$. Therefore,  assumption {\bf (H)} is satisfied.

\paragraph{Example 5.2 (Diffusions on Heisenberg group).} Consider the following vector fields on $\R^3$:
$$U_1(x)=  \pp_1 - \ff {x_2} 2  \pp_3,\ \  U_2(x)= \pp_2 +\ff {x_1} 2 \pp_3,\ \ x=(x_1,x_2,x_3)\in\R^3.$$ Then the H\"ormander condition is satisfied. The Markov semigroup $\e^{t\DD_0}$ generated by  the Kohn-Laplacian $\DD_0:=U_1^2 +U_2^2$   has a strictly positive smooth density $p_t$:
$$(\e^{t\DD_0}f)(0)= \int_{\R^3} p_t(x) f(x)\d x,\ \ f\in \B_b(\R^3), t>0.$$
For fixed $t_0>0$, let  $V(x)=-\log p_{t_0}(x)$ so that $\mu_0(\d x):= \e^{-V(x)}\d x= p_{t_0}(x)\d x.$   Obviously, $1\in H_\si^{1,2}(\mu_0)$ with $\EE_0(1,1)=0$. Moreover, according to \cite[Corollary 1.2]{Li}, there exists a universal constant $c>0$ such that the log-Sobolev inequality \eqref{LS} holds for $\kk=t_0c$ and $\bb=0$ (see \cite{W*} for more results on functional inequalities).
So, assumption {\bf (H)} holds.

\

In the next two subsections, we   prove Theorem \ref{T3.1} for $\tau=0$ and $\tau>0$ respectively.

\subsection{Proof of Theorem \ref{T3.1} for $\tau=0$}

When $\tau=0$, Theorem \ref{T3.1} reduces to the following result where  $S_t^Z$ is replaced by $P_t^Z$ for   notation consistency with  Theorem \ref{T1.1}.

\beg{thm}\label{T3.2} Assume {\bf (H)} and let $\tau=0$. For any $\ll>\kk$ and $p\in \big(1, \ff{\ss\ll}{\ss\kk}\big)$, there exists a constant $C_{p,\ll}$ such that
 $\mu_0(\e^{\ll |Z|^2})<\infty $  implies that $\scr P_Z$ has a unique element   $\nu=\rr\mu_0$, where $\rr$ is continuous, strictly positive satisfying
 $\log \rr, \rr^{p/2} \in H_\si^{1,2}(\mu_0)$,   $\eqref{SI2}$ and
\beq\label{DD} \mu_0\big(|\si^*\nn \rr^{\ff p 2}|^2+ \rr^p\big)\le C_{p,\ll} \{\mu_0(\e^{\ll |Z|^2})\}^{C_{p,\ll}}.\end{equation}\end{thm}

We first prove this result for bounded $Z$ with compact support. Let $\scr L_Z= \scr L_0+(\si Z)\cdot\nn,$ where $\scr L_0$ is in \eqref{SL}.
Then  an  invariant probability measure $\nu$ of $P_t^Z$ solves the equation  $\scr L_Z^*\nu=0$ in the sense that
$$\int_{\R^d} \scr L_Zf\d\nu=0, \ \ f\in C_0^\infty(\R^d).$$

The following lemma  extends  \cite[Lemma 4.2]{W16} to the present degenerate case.

\beg{lem}\label{L00} Assume  {\bf (H)}  and let    $\tau=0$. If   $|Z|\in \cap_{p\in [1,\infty)} L^p_{loc}(\R^d)$ and  $\nu:=\rr\mu_0$ is a probability measure such that   $L_Z^*\nu=0$, then  $\rr$ has a continuous, strictly positive version.
 If moreover $Z$ is bounded and has compact support, then $\rr\in H_\si^{1,2}(\mu_0)$ and
\beq\label{FF}  \int_{\R^d} \<\si^*\nn f, \si^*\nn\rr\>\d\mu_0 =  2\int_{\R^d} \<Z, \si^*\nn f\>\rr\d\mu_0,\ \ f\in H_\si^{1,2}(\mu_0).\end{equation}
\end{lem}

\beg{proof}


We first  prove that $\rr$ has a continuous, strictly positive version  using results in \cite{NN}. Let $v_i$ stand for the $i$-th component of a vector $v$, and let
$$A_i(x,u,\xi)= \xi_i -2 uZ_i(x),\ \ x\in\R^d, u\in \R, \xi\in \R^m,\ \ 1\le i\le m.$$
It follows from  the integration by parts formula and $\scr L_Z^*\nu=0$ that
\beq\label{*K}\beg{split}  &\int_{\R^d} \sum_{i=1}^m A_i(\cdot, \rr, (U_1\rr,\cdots, U_m\rr)) U_i f \d\mu_0= \sum_{i=1}^m\int_{\R^d} (U_i\rr-2 \rr Z_i) U_i f \d\mu_0\\
&=:-\sum_{i=1}^m\int_{\R^d} \Big({\rm div}\big\{(\e^{-V}U_i f)U_i\} + 2\e^{-V} Z_i U_i f\Big)\rr \d x= -2\nu(L_Zf)=0,\ \ f\in C_0^\infty(\R^d).\end{split}\end{equation}
Obviously,
\beg{align*} &|A(x,u,\xi)|:=\sum_{i=1}^m |A_i(x,u,\xi)|\le 2 |u| \cdot|Z|(x) +|\xi|,\\
&A(x,u,\xi)\cdot \xi\ge |\xi|^2 -|Z|(x)|u|\cdot|\xi|\ge \ff 1 2  |\xi|^2 -2 |u|^2|Z|^2(x),  \end{align*}
where $|Z|\in L_{loc}^p(\mu)$ for any $p\in [1,\infty)$.
Then by \cite[Theorem 3.1 and Theorem 3.35]{NN}, $\rr$ has a locally H\"older continuous version (denoted again by $\rr$) with respect to the intrinsic distance induced by $\si$. By the H\"ormander condition, the intrinsic distance induces the classical topology in $\R^d$, so that this version $\rr$ is continuous. Moreover,  for any compact set $K$ there exists a constant $C(K)>0$ such that
$$\sup_K \rr\le C(K)+C(K)\inf_K \rr.$$ Since the equation \eqref{*K} is linear in $\rr$, this inequality also holds for  $n\rr$ replacing $\rr$,  so that
$$\sup_K \rr\le \ff { C(K)}n+C(K)\inf_K \rr,\ \ n\ge 1.$$ Letting $n\to\infty$ we obtain
$$\sup_K \rr\le  C(K)\inf_K \rr,$$ which implies that $\rr$ is strictly positive since $\mu_0(\rr)=1$.

Now, let $Z$ be bounded with   compact support. Since $\rr$ is locally bounded due to the continuity,   $\rr Z$ is globally bounded. In particular, $\rr Z\in L^2(\R^d\to \R^d;\mu_0)$.  By \eqref{P} and the completion of $H_\si^{1,2}(\mu_0)$, $\{\si^*\nn h:\ h\in H_\si^{1,2}(\mu_0)\}$ is a closed subspace of $L^2(\R^d\to \R^d;\mu_0)$. Let $h\in H_\si^{1,2}(\mu_0)$ such that $\si^*\nn h$ is the orthogonal projection of $\rr Z$ on this subspace.
Since $\scr L_Z^*\nu=0$, we have
\beg{align*} &\mu_0(\rr \scr L_0f)=   \nu(\scr L_Z f)-\nu(\<Z,\si^*\nn f\>)\\
&= -\mu_0(\<\rr Z,\si^*\nn f\>)=-\EE_0(h,f)=\mu_0(h\scr L_0 f), \ \ f\in C_0^\infty(\R^d).\end{align*}
By \eqref{P}, $\{\scr L_0 f: f\in C_0^\infty(\R^d)\}$ is dense in $\{f\in L^2(\mu_0): \mu_0(f)=0\}$, so this   implies $\rr =h+c$ for some constant $c$. Hence,    $\rr\in H_\si^{1,2}(\mu_0)$ and
$$\ff 1 2 \mu_0(\<\si^*\nn \rr, \si^*\nn f\>)=    \mu_0(\rr \<Z,\si^*\nn f\>),\ \ f\in C_0^\infty(\R^d).$$ Due to the boundedness of $\rr Z$, this is equivalent to   \eqref{FF}.
\end{proof}

\beg{lem}\label{L3.3'} Assume {\bf (H)} and let $\tau=0$.  For any $\ll>\kk$ and $p\in \big(1, \ff{\ss\ll}{\ss\kk}\big)$, there exists a constant $C_{p,\ll}$ such that  for any bounded $Z$ with compact support, if a probability measure $\nu:=\rr\mu_0$ solves $L_Z^*\nu=0$, then  $\rr$  has a continuous,  strictly positive version such that  $\log \rr, \rr^{p/2} \in H_\si^{1,2}(\mu_0)$,
  and $\eqref{SI2}$, $\eqref{DD}$ hold.
 \end{lem}

 \beg{proof}  By  Lemma \ref{L00}, $\rr$ has a continuous, strictly positive version such that $\rr\in H^{1,2}_\si(\mu_0)$ and \eqref{FF} holds.   According to step (a) in the proof of \cite[Theorem 2.3]{W16},  \eqref{FF} implies   $\log\rr, \ss\rr \in H_\si^{1,2}(\mu_0)$ and \eqref{SI2}. To prove \eqref{DD},    let $f_n:= (n^{-1}+ \rr\land n)^{p-1},\ n\ge 1. $ Then $f_n\in H_\si^{1,2}(\mu_0)$ and by \eqref{FF},
\beq\label{D1}\beg{split} I_n&:= \mu_0\big( \big|\si^* \nn (\rr\land n+n^{-1})^{\ff{p}2}\big|^2\big)
 = \ff{p^2}{4(p-1)} \int_{\R^d}  \<\si^*\nn f_n, \si^*\nn\rr\>\d\mu_0  \\
 &=\ff{p^2}{2(p-1)} \int_{\R^d}  \<Z, \si^*\nn f_n\> \rr \d\mu_0
 \le \ff{p^2}{2(p-1)}\mu_0\big(\{\rr\land n\}|Z|\cdot |\si^* \nn f_n| \big) \\
&\le   p \mu_0\big( |Z|(\rr\land n+n^{-1})^{\ff {p}2}\big |\si^*\nn (\rr\land n+n^{-1})^{\ff p2}\big|\big)\\
&\le   p  \ss{I_n \mu_0\big(|Z|^2 (\rr\land n+n^{-1})^p\big)}.\end{split}\end{equation}
Next,  it follows from \eqref{Young}  and \eqref{LS} that
\beq\label{D2} \beg{split}&\mu_0\big(|Z|^2 (\rr\land n+n^{-1})^p\big)-\ff {\mu_0((\rr\land n+n^{-1})^p)} \ll \log \mu_0(\e^{\ll|Z|^2})\\
&\le \ff 1 \ll \mu_0\bigg((\rr\land n+n^{-1})^p\log \ff{(\rr\land n+n^{-1})^p}{\mu_0((\rr\land n+n^{-1})^p)}\bigg)\\
&\le \ff \kk \ll I_n +\ff{\bb}\ll \mu_0\big((\rr\land n+n^{-1})^p\big).\end{split}\end{equation}
Noting that the log-Sobolev inequality \eqref{LS} implies the super Poincar\'e inequality (see \cite{W00a} or \cite{W00b})
\beq\label{SP} \mu_0(f^2)\le r\mu_0(|\si^*\nn f|^2)+\e^{c(1+ r^{-1})}\mu_0(|f|)^2,\ \ r>0, f\in H_\si^{1,2}(\mu_0)\end{equation}
for some constant $c>0$, we obtain
\beq\label{D3} \mu_0\big((\rr\land n+n^{-1})^p\big)\le r I_n +\e^{c(1+r^{-1})} \mu_0\big((\rr\land n+n^{-1})^{\ff {p}2}\big)^2,\ \ r>0.\end{equation}
Combining \eqref{D1}-\eqref{D3}, we arrive at
\beq\label{D4}\beg{split}  & I_n  \le p^2\mu_0(|Z|^2(\rr\land n+n^{-1})^p)\\
&\le p^2  \Big(\ff \kk \ll I_n + \ff{\bb+\log \mu_0(\e^{\ll|Z|^2})}\ll \mu_0\big((\rr\land n+n^{-1})^p\big)\Big)\\
&\le  p^2  \Big(\ff \kk \ll I_n + \ff{\bb+\log \mu_0(\e^{\ll|Z|^2})}\ll \Big\{rI_n +\e^{c(1+r^{-1})}  \mu_0\big((\rr\land n+n^{-1})^{\ff{p}2}\big)^2\Big\}\Big)\\
& = \ff{p^2}{\ll} \Big[\big\{\kk+r(\bb+\log \mu_0(\e^{\ll|Z|^2}))\}  I_n\\
&\qquad\qquad\qquad  + \big\{\e^{c(1+r^{-1})}(\bb+\log \mu_0(\e^{\ll|Z|^2})) \big\}  \mu_0\big((\rr\land n+n^{-1})^{\ff{p}2}\big)^2 \Big],\ \ r>0.
\end{split}\end{equation}
Now, we are ready to complete the proof by induction in $p$ as follows.

(i) Assume that $p\le 2$. Since $\mu_0(\rr)=1$, we have
\beq\label{DD3} \mu_0\big((\rr\land n+n^{-1})^{\ff {p}2}\big)^2\le \{\mu_0(\rr+1)\}^{p}\le 2^2=4,\ \ n\ge 1.\end{equation}
So, \eqref{D4} implies
$$I_n\le \ff{p^2}{\ll} \Big[\big\{\kk+r(\bb+\log \mu_0(\e^{\ll|Z|^2}))\}  I_n
  + 4\big\{\e^{c(1+r^{-1})}(\bb+\log \mu_0(\e^{\ll|Z|^2})) \big\}   \Big],\ \ r>0.$$ Since $p^2 <  \ff{\ll}\kk,$ letting
  $$r= r_{p,\ll}:= \ff{\ll -p^2 \kk}{2p^2 (\bb+\log \mu_0(\e^{\ll|Z|^2}))}>0, $$
we obtain
  $$I_n\le \ff{8 p^2\e^{c(1+r_{p,\ll})}(\bb+\log \mu_0(\e^{\ll|Z|^2}))  }{\ll-p^2\kk}\le \aa_{p,\ll}\mu_0(\e^{\ll|Z|^2})^{\aa_{p,\ll}},\ \ \ n\ge 1$$ for some constant $\aa_{p,\ll}>0$.
Combining this with \eqref{D3} for e.g. $r=1$ and \eqref{DD3}, we may find out a constant $C_{p,\ll}>0$ such that
\beq\label{D5} I_n+\mu((\rr\land n)^p)\le C_{p,\ll}\mu_0(\e^{\ll|Z|^2})^{C_{p,\ll}},\ \ n\ge 1.\end{equation}
Letting $n\to\infty$ we conclude that   $\rr^{p/2}\in H_\si^{1,2}(\mu_0)$ and  \eqref{DD} holds.

(ii) Assume that \eqref{DD} holds with $p\in (1, \ff{\ss\ll}{\ss\kk})\cap [0,k]$ for some $k\ge 2$, we aim to prove it for $p\in (1, \ff{\ss\ll}{\ss\kk})\cap [0,2k]$. It suffices to consider   $p\in (1, \ff{\ss\ll}{\ss\kk})\cap (k,2k]$ and $k<\ff{\ss\ll}{\ss\kk}$. In this case,  by the assumption  there exists a constant $\aa_{k,\ll} \ge 1$ such that
$$\mu_0(\rr^k)\le \aa_{k,\ll}  \mu_0(\e^{\ll|Z|^2})^{\aa_{k,\ll}}.$$ Since $p\le 2k$, we have
$$\mu_0((\rr\land n+n^{-1})^{p/2})^2\le \{\mu_0(\rr^k)\}^{p/k}\le \aa_{k,\ll}^{p}  \mu_0(\e^{\ll|Z|^2})^{p\aa_{k,\ll}}.$$
Substituting this into  \eqref{D4} and repeating the argument in (i), we prove \eqref{DD} for some constant $C_{p,\ll}>0.$
\end{proof}

\beg{proof}[Proof of Theorem \ref{T3.2}] By the H\"ormander theorem, {\bf (H)} implies that $P_t^0$ has a strictly positive density with respect to $\mu_0$. So, the uniqueness of $\nu\in \scr P_Z$ follows from Theorem \ref{T0}(1).

To prove the existence,
for any $n\ge 1$, let $Z_n= Z1_{\{|Z|+|\cdot|\le n\}}.$ Then $Z_n$ is bounded with compact support. By Theorem \ref{T2.2} and Lemma \ref{L3.3'}, $P_t^{Z_n}$ has an invariant probability measure $\nu_n=\rr_n\mu_0$, where $\rr_n$ is continuous, strictly positive  such that   $\log\rr_n, \rr_n^{p/2}\in H_\si^{1,2}(\mu_0)$ and \eqref{SI2}, \eqref{DD} hold for $\rr_n$ replacing $\rr.$ In particular,    $\{\rr_n^{p/2}\}_{n\ge 1}$ is bounded  in $H_\si^{1,2}(\mu_0)$. Then, as explained in step (b) in the proof of \cite[Theorem 2.3(1)]{W16},  the defective log-Sobolev inequality \eqref{LS} implies that  $\{\rr_n^{p/2}\}_{n\ge 1}$ is relatively compact in $L^2(\mu_0)$. So,    up to a subsequence,  $\rr_n^{p/2}\to \rr^{p/2}$ in $L^2(\mu_0)$ for some probability density $\rr$ with respect to $\mu_0$. Moreover, $\log\rr, \rr^{p/2}\in H_\si^{1,2}(\mu_0)$ and \eqref{SI2}, \eqref{DD} hold.
It remains to prove that $\nu:=\rr\mu_0$ is an invariant probability measure of $P_t^Z$, where $\rr$ has a continuous, strictly positive version according to Lemma \ref{L00}.

For any $f\in \B_b(\R^d)$, by $\nu_n(P_t^{Z_n}f)=\nu_n(f)$, $\rr_n\to \rr$ in $L^1(\mu_0)$ and the boundedness of $\mu_0(\rr_n^p)$, we obtain
\beq\label{D5}\beg{split}& \big|\nu(P_t^Z f)-\nu(f)\big| = \Big|\lim_{n\to\infty} \mu_0(\rr_n (P_t^Z f-f))\Big|\\
&\le  \limsup_{n\to\infty}  \mu_0\big( \rr_n |P_t^Z f-P_t^{Z_n}f|\big)\le C \limsup_{n\to\infty}  \mu_0\big(  |P_t^Z f-P_t^{Z_n}f|^{q}\big)^{\ff{1}q} \end{split}\end{equation} for some constant $C>0$ and $q:=\ff p{p-1}$. By \eqref{SM} and   $\tau=0$,    we have
\beg{equation}\label{D6}\beg{split} &\mu_0\big(  |P_t^Z f-P_t^{Z_n}f|^{q}\big) \\
&= \int_{\R^d} \big|\E \big[f(X^x(t))\big\{\e^{\int_0^t\<Z(X^x(s)), \d W(s)\>-
 \ff 1 2\int_0^t |Z(X^x(s))|^2\d s } \\
 &\qquad\qquad - \e^{\int_0^t\<Z_n(X^x(s)), \d W(s)\>-
 \ff 1 2\int_0^t |Z_n(X^x(s))|^2\d s } \big\}\big]\big|^q\mu_0(\d x)\\
 &\le \E \big[\big|f(X^{\mu_0}(t))\big\{\e^{\int_0^t\<Z(X^{\mu_0}(s)), \d W(s)\>-
 \ff 1 2\int_0^t |Z(X^{\mu_0}(s))|^2\d s } \\
 &\qquad \qquad- \e^{\int_0^t\<Z_n(X^{\mu_0}(s)), \d W(s)\>-
 \ff 1 2\int_0^t |Z_n(X^{\mu_0}(s))|^2\d s } \big\}\big|^q\big]\\
 &\le \|f\|_\infty^q \big(\E \e^{2q\int_0^t\<Z(X^{\mu_0}(s)), \d W(s)\>-
 q\int_0^t |Z(X^{\mu_0}(s))|^2\d s }\big)^{\ff 1 2} \\
 &\qquad \qquad\times \Big(\E \big|\e^{\int_0^t\<(Z-Z_n)(X^{\mu_0}(s)), \d W(s)\>-
\ff 1 2  \int_0^t |(Z-Z_n)(X^{\mu_0}(s))|^2\d s }-1\big|^{2q}\Big)^{\ff 1 2}. \end{split} \end{equation} Since $\mu_0(\e^{\ll|Z|^2})<\infty$, for any $\aa>1$ and measurable function $g$ with $|g|\le |Z|$,
\beg{align*}& \E \e^{2\aa\int_0^t\<g(X^{\mu_0}(s)), \d W(s)\>-\aa \int_0^t |g(X^{\mu_0}(s))|^2\d s} \\
&\le \Big\{\big(\E\e^{4\aa \int_0^t\<g(X^{\mu_0}(s)), \d W(s)\>-8\aa^2 \int_0^t |g(X^{\mu_0}(s))|^2\d s} \big) \E\e^{(8\aa^2-2\aa)\int_0^t |g(X^{\mu_0}(s))|^2\d s}\Big\}^{\ff 1 2}\\
&\le \big(\mu_0(\e^{(8\aa^2-2\aa)t|Z|^2})\big)^{\ff 1 2}\le \ss{\mu_0(\e^{\ll|Z|^2})}<\infty,\ \ t\le \ff\ll{8\aa^2-2\aa}.\end{align*}
So, for small enough $t>0$, by \eqref{D6} and  the  dominated convergence theorem we prove
$\lim_{n\to\infty} \mu_0\big(|P_t^Z f-P_t^{Z_n}f|^{q}\big)=0.$ Then    \eqref{D5} implies $\nu(P_t^Zf)=\nu(f)$ for small $t>0$ and all $f\in \B_b(\R^d)$. Therefore, $\nu$ is an invariant probability measure of $P_t^Z$.
\end{proof}

\subsection{Proof of Theorem \ref{T3.1} for  $\tau>0$}

 Again we start from   bounded $Z$.

\beg{lem}\label{LN} Assume {\bf (H)} and let $Z$ be bounded. Then  $S_t^Z$ has an invariant probability measure $\nu=\rr\mu$ such that $\rr\in \cap_{p\in (1,\infty)} L^p(\mu)$. Moreover,
    $ \nu_\theta=\nu_0$ for $\theta\in [-\tau, 0]$,  and it has a continuous, strictly positive density $\rr_0$ with respect to $\mu_0$
such that $\log\rr_0, \rr_0^{p/2}\in H_\si^{1,2}(\mu_0)$ for any $p>1$. Moreover, for any $\ll>\kk$ and $p\in (1, \ff {\ss\ll}{\ss\kk}),$  there exists $k=k(p,\ll)>0$  such that $\eqref{SI3}$ holds.   \end{lem}

\beg{proof}     Since \eqref{LS} implies the hyperboundedness of $P_t^0$, the existence   of invariant probability measure $\nu=\rr\mu$ is ensured by Theorem  \ref{T2.2}, which is the weak limit of a subsequence of $\{\nu_n\}_{n\ge 1}$ in \eqref{SEQ}.  Below we first prove the assertion on the marginal density, then prove $\rr\in L^p(\mu)$ for all $p>1$.

(a) The marginal density. Let $\theta\in [-\tau,0]$. Since  $\nu$ is $S_t^Z$-invariant,  for any $f\in \B_b(\R^d)$ and
$  f_\theta(\xi):= f(\xi(\theta))$, we have
\beg{align*} &\nu_\theta (f)=   \nu(f_\theta) = \nu(S_{-\theta}^Z  f_\theta) = \int_\C\E\big[f_\theta(X_{-\theta}^\xi) R^\xi(-\theta)\big]\nu(\d\xi)\\
&= \int_\C f(X^\xi(-\theta+\theta))\nu(\d \xi)=\int_\C f(\xi(0))\nu(\d\xi)
=\nu_0(f).\end{align*} Therefore, $\nu_\theta=\nu_0$.

Moreover,
for any $f\in C_0^\infty(\R^d)$ and $ f_0(\xi):= f(\xi(0))$, we have
$$  \E \big[f(X^{\nu}(t))R^{\nu}(t)\big] =\nu(S_t^Z   f_0)=\nu(f_0)= \nu_0 (f),\ \ t\ge 0.$$ On the other hand, let   $Z|_{\xi(0)}= \nu(Z|\xi(0))$ be the regular conditional expectation of $Z$ under probability $\nu$ given $\xi(0)$.   By It\^o's formula,
 $$ \E \big[f(X^{\nu}(t))R^{\nu}(t)\big]-\nu_0(f) = \int_0^t \E \big[( L_{Z|_{\xi(0)} } f)(X^{\nu}(s))R^{\nu}(s)\big]\d s
 =\int_0^t \nu_0(L_{Z|_{\xi(0)}}f)\d s.$$ Therefore,
 $\nu_0(L_{Z|_{\xi(0)}}f)=0$  for all $f\in C_0^\infty(\R^d)$, i.e. $L_{Z|_{\xi(0)}}^*\nu_0=0$. We then finish the proof by considering the following two situations.

 (i) $Z$ is supported on a bounded subset of $\C$. Then $Z|_{\xi(0)}$ has compact support.
Since  by Jensen's inequality
 $$\mu_0( \e^{\ll|Z|_{\xi(0)}|^2})= \mu\big(\e^{\ll|\mu(Z|\xi(0))|^2}\big)\le \mu\big[\mu(\e^{\ll|Z|^2}|\xi(0))\big]= \mu( \e^{\ll|Z|^2}) <\infty,$$    the desired assertion on $\rr_0$  follows from   Lemma  \ref{L3.3'}.

 (ii) In general, let $Z^{\<n\>}= Z1_{\{|\cdot|\le n\}}$ for $n\ge 1$. Then for every $n\ge 1$, $S_t^{Z^{\<n\>}}$ has an invariant probability measure $\nu^{\<n\>}= \rr^{\<n\>}\mu$ with the marginal density $\rr_0^{\<n\>}$ satisfying \eqref{SI3} in place of $\rr$. As shown in the proof of Theorem \ref{T2.2} that up to a subsequence $\rr^{\<n\>}\to\rr$ weakly in $L^1(\mu)$,  $\nu:=\rr \mu$ is an invariant probability measure of $S_t^Z$, and $\rr_0$  satisfies  \eqref{SI3}. Note that by Lemma \ref{L00}, $\rr_0$ has a continuous, strictly positive version, so that
 the Poincar\'e inequality \eqref{P}   implies $\log \rr_0\in L^2(\mu_0)$, see step (a) in the proof of  \cite[Theorem 2.3(1)]{W16} for details.
 Therefore, after proving $\rr\in L^p(\mu)$ (hence, $\rr_0\in L^p(\mu_0)$) in the next step, we conclude that $\log \rr_0, \rr_0^{\ff p 2}\in H_\si^{1,2}(\mu_0)$.

 (b) $\rr\in L^p(\mu)$ for $p>1$. Let $f\ge 0$ with $\mu(f^{\ff p{p-1}})\le 1$. Since $Z$ is bounded, for any $\aa>1$ there exists a constant $c(\aa)>0$ such that $\E[ (R^\nu(\tau))^\aa]\le c(\aa)<\infty$. Let $q\in (1,\ff p{p-1})$. Combining this with \eqref{3*} and \eqref{SM} and using $\rr_0\in L^\aa(\mu_0)$ for any $\aa>1$, we obtain
\beq\label{HH}\beg{split} \nu(f)&= \nu( S_\tau^Zf) = \E [f(X_\tau^\nu)R^\nu(\tau)] \le c_1\big\{\E [f^q(X_\tau^\nu)]\big\}^{\ff 1 q}\\
&= c_1 \big\{\mu_0(\rr_0 S_\tau f^q)\big\}^{\ff 1 q} \le c_1\big\{\mu_0\big((S_\tau f^q)^{\ff p{q(p-1)}}\big)\big\}^{\ff {p-1}p} \big\{\mu_0(\rr_0^{\ff{p}{p-(p-1)q}})\big\}^{\ff{p-(p-1)q}{pq}}\\
&\le c_2 \big\{\mu(f^{\ff p{p-1}})\big\}^{\ff {p-1}p},\ \ f\in \B^+(\C)\end{split} \end{equation} for some constants $c_1,c_2>0$. Therefore, $\rr\in L^p(\mu).$
 \end{proof}

\beg{proof}[Proof of Theorem \ref{T3.1}] By H\"ormander's theorem, {\bf (H)} implies that $P_t^0$ has a strictly positive density with respect to $\mu_0$. So, the uniqueness of $\nu$ follows from Theorem \ref{T0}(1). Below, we prove the existence and assertions (1) and (2).

(a) We first assume that $Z$ is bounded. By Lemma \ref{LN}, it remains to prove \eqref{SI4} for $p<p_\ll:=\ff{q_\ll}{q_\ll-1}$ but close enough to $p_\ll$. Since $\ll>\ll_{\kk,\tau}$, we have
$$p_3:= \ff 1 4 \big(1+\ss{1+8\ll/\tau}\big) >\ff{4\ss\ll}{\ss\ll-\ss\kk}$$ and
$$q_\ll= \ff 1 2 \bigg(p_3 -\ss{p_3^2 -\ff{4p_3\ss\ll }{\ss\ll-\ss\kk}}\bigg)<\ff {p_3}2.$$
Since $\ff{p_\ll}{p_\ll-1}= q_\ll$, when $p\in (1,p_\ll)$ is close enough to $p_\ll$  we have
$$q_\ll <\ff p{p-1} <\ff 1 2 \bigg(p_3 +\ss{p_3^2 -\ff{4p_3\ss\ll }{\ss\ll-\ss\kk}}\bigg)<p_3,$$ so that
\beq\label{FS0} \Big(\ff p {p-1}\Big)^2-\ff{p_3p}{p-1} + \ff{p_3\ss\ll }{\ss\ll-\ss\kk}<0\end{equation} and
$$p_1:= p_3 (p-1)>p,\ \ p_2:= \ff{p_1}{p_1-p}=\ff{p_3(p-1)}{p_3(p-1)-p}>1.$$
It is easy to see that \eqref{FS0} is equivalent to
$$\ff{\ff{\ss\ll}{\ss\kk} -1}{p_2\ff{\ss\ll}{\ss\kk}}= \ff{\ss\ll-\ss\kk}{\ss\ll} \Big(1-\ff p{p_3(p-1)}\Big) >\ff {p-1}p.$$
Then  there exists $\theta\in (1,\ff{\ss\ll}{\ss\kk})$ depending on $p,\ll$ such that $\ff{\theta -1}{p_2\theta }=\ff {p-1}p.$
By \eqref{3*}  and   Lemma \ref{LN},  there exists a constant $k=k(p,\ll)$ such that
\beq\label{H1} \beg{split}&\int_\C  \rr(\xi) \E[f^{p_2}(X_\tau^\xi)] \mu(\d \xi)=\mu_0(\rr_0 S_\tau f^{p_2})
 \le \{\mu_0(\rr_0^\theta)\}^{\ff 1 \theta} \{\mu_0(S_\tau f^{\ff{p_2\theta}{\theta-1}})\}^{\ff{\theta-1}\theta}\\
 &\le \{k\mu(\e^{\ll|Z|^2})^k\}^{p_2} \{\mu(f^{\ff{p_2\theta}{\theta-1}})\}^{\ff{\theta-1}\theta}\le \{k\mu(\e^{\ll|Z|^2})^k\}^{p_2},\ \ f\ge 0, \mu(f^{\ff p{p-1}})\le 1.\end{split}\end{equation}
Noting that
$$\ff 1 {p_1}+\ff 1 {p_2}+\ff 1 {p_3}= \ff{1+p_3(p-1)-p+p-1}{p_3(p-1)}=1,$$
by   H\"older's inequality and \eqref{H1}, we obtain
\beq\label{H0} \beg{split}\nu(f)&=\nu(P_\tau^Zf) =\int_\C  \rr(\xi) \E[f(X_\tau^\xi)R^\xi(\tau)] \mu(\d \xi)\\
&\le \mu(\rr^p)^{\ff 1 {p_1}} \bigg(\int_\C\rr(\xi) \E[f^{p_2}(X_\tau^\xi)] \mu(\d \xi)\bigg)^{\ff 1 {p_2}} \big(\E[(R^{\mu}(\tau))^{p_3}]\big)^{\ff 1 {p_3}}\\
&\le k\mu(\e^{\ll|Z|^2})^k \mu(\rr^p)^{\ff 1 {p_1}}\big(\E[(R^{\mu}(\tau))^{p_3}]\big)^{\ff 1 {p_3}},\ \ f\ge 0, \mu(f^{\ff p{p-1}})\le 1.\end{split}\end{equation}
To estimate $\E [(R^\mu(\tau))^{p_3}]$, recall that for any continuous martingale $M(t)$   we have
$$\E \e^{p_3M(\tau)-\ff{p_3}2\<M\>(\tau)} \le \big(\E\e^{2 p_3 M(\tau) -  2p_3^2  \<M\>(\tau)}\big)^{\ff 1 2} \big(\E \e^{ p_3(2 p_3-1)\<M\>(\tau)}\big)^{\ff 1 2}
 \le \big(\E \e^{  p_3(2 p_3-1) \<M\>(\tau)}\big)^{\ff 1 2}.$$
Taking  $M(t)=\int_0^t\<Z(X_s^\mu), \d W(s)\>$ and  noting that $p_3(2p_3-1)\tau=\ll$  by the definition of $p_3$, we obtain
\beg{align*} &\big(\E(R^\mu(\tau))^{p_3}\big)^2 \le \E\big[\e^{ p_3 (2p_3-1)  \int_0^\tau |Z(X_s^\mu)|^2\d s}\big]
\\
&\le \ff 1\tau \int_0^\tau \E\big[\e^{ p_3 (2p_3-1) \tau   |Z(X_s^\mu)|^2 }\big]\d s  =\mu(\e^{\ll|Z|^2}).\end{align*}
  Combining this with \eqref{H0},  we arrive at
$$\mu(\rr f)=\nu(f)\le \mu(\rr^p)^{\ff 1 {p_1}} k\mu(\e^{\ll |Z|^2})^k,\ \ f\ge 0, \mu(f^{\ff p{p-1}})\le 1$$ for some constant $k=k(p,\ll)$. Since $p_1>p$ and $\rr\in L^p(\mu)$ due to Lemma \ref{LN}, this implies the desired estimate \eqref{SI4}.

(b)  In general, for any $n\ge 1$, let $Z^{\<n\>}=Z1_{\{|Z|\le n\}}.$ Then $S_t^{Z^{\<n\>}}$ has an invariant probability measure $\nu^{\<n\>}:=\rr^{\<n\>}\mu$ such that \eqref{SI4} and \eqref{SI3} hold for $\rr^{\<n\>}$ and $\rr_0^{\<n\>}$ replacing
$\rr$ and $\rr_0$ with constants independent of $n$. In particular, $\rr^{\<n\>}$ converges weakly in $L^1(\mu)$ to some $\rr$, and as shown in the proof of Theorem \ref{T2.2} that $\nu:= \rr\mu$ is an invariant probability measure of $S_t^Z$ satisfying \eqref{SI4} and \eqref{SI3}.
Moreover, applying Lemma \ref{L00} to the marginal distribution $\nu_0$ (recall that $L^*_{Z|_{\xi(0)}}\nu_0=0$), we conclude that $\rr_0$ has a continuous, strictly positive version.

\end{proof}

\subsection{The infinite-dimensional case}

By finite-dimensional approximations, it is easy to extend Theorem \ref{T3.1} to the infinite-dimensional case. For simplicity, here we only consider an Ornstein-Uhlenbeck type reference process on $\H$.

 Let $W(t)$ be the cylindrical Brownian motion on $\H$, and let $L, \si$ be self-adjoint operators such that for some orthonormal basis $\{e_i\}_{i\ge 1}$ of $\H$
$$Le_i =\ll_i e_i, \ \ \ \si e_i= q_ie_i,\ \ i\ge 1$$ holds for some constants $\ll_i,q_i$ satisfying
$$\ll_1=\inf_{i\ge 1}\ll_i >0,\ \ \inf_{i\ge 1} q_i^2 >0, \ \ \sum_{i\ge 1} \ff{q_i^2}{\ll_i^\dd}<\infty \ \text{for\ some}\ \dd\in (0,1).$$ Then for any initial point,
the SDE
\beq\label{SLE}\d X(t)= -L X(t) \d t +\si \d W(t) \end{equation} has a unique continuous mild solution, and the associated Markov semigroup $P_t^0$ is symmetric in $L^2(\mu_0)$ for $\mu_0$ being the centered Gauss  measure on $\H$ of covariance operator $Q$ with $Qe_i:= \ff {q_i^2}{2\ll_i}e_i,\ i\ge 1.$  When $\tau>0,$ let $\mu$ be the distribution of $X_\tau^{\mu_0}$ as introduced in Section 1.

Next,  according to \cite{Gross}, we have the following log-Sobolev inequality
\beq\label{LS*} \mu_0(f^2\log f^2)\le \ff 1 {\ll_1} \mu_0(|\si^*\nn f|^2),\ \ f\in \F C_0^\infty, \mu_0(f^2)=1,\end{equation}
where $\F C_0^\infty:= \{x\mapsto f(\<x,e_1\>, \cdots, \<x, e_n\>):\ n\ge 1, f\in C_0^\infty(\R^n)\}$  is the class of smooth cylindrical functions.

Below we   extend Theorem \ref{T3.1} to the SDE
\beq\label{SDE'} \d X(t) = \big\{\si Z(X_t) -L X(t)\big\}\d t +\si \d W(t)\end{equation}  by using finite-dimensional approximations.

 For any $n\ge 1$, let $\pi_n:\H\to \H_n:={\rm span}\{e_1,\cdots, e_n\}$ be the orthogonal projection. Then $\pi_n X(t)$ is a Markov process on $\H_n$ which is symmetric with respect to $\mu_0^{(n)}:= \mu_0\circ \pi_n^{-1}$, and
 \eqref{LS*} implies the same log-Sobolev inequality for $\mu_0^{(n)}$ on $\H_n$ replacing $\mu_0$ on $\H$. Let $Z: \C\to \H$ be measurable satisfying conditions in Theorem \ref{T3.1}. Then  $Z_n =Z|_{\C_n}$ also satisfies these conditions, where $\C_n:=C([-\tau,0]; \H_n)\subset \C$. So,  letting $\mu^{(n)} $ be the marginal distribution of $\mu$ on $\C_n$, the  corresponding finite-dimensional Markov semigroup $S_t^{Z_n}$   has an invariant probability measure $\nu^{(n)}=\rr^{(n)}\mu^{(n)}$ with $\rr^{(n)} $ and $\rr^{(n)}_0$ satisfying \eqref{SI4} and  \eqref{SI3} respectively. Thus, up to a subsequence, $\rr^{(n)}\circ\pi_n\to \rr$ weakly in $L^p(\mu)$, and
$(S_t^{Z_n} f)\circ \pi_n \to S_t^Z f$ in $L^{\ff p{p-1}}(\mu)$ for any bounded cylindrical function  $f$ on $\H$.   Therefore, $\nu:=\rr\mu\in \scr P_Z$ with $\rr$ and $\rr_0$  satisfying \eqref{SI4} and \eqref{SI3} respectively.

Moreover, let $P_t^0$ be the   Markov semigroup of the linear equation \eqref{SLE}. According to e.g. \cite[Theorem 3.2.1]{Wbook}, $P_t^0$ satisfies the following Harnack inequality for some constant $C>0$:
$$(P_t^0f(x))^p\le (P_t^0 f^p(y))\exp\bigg[\ff{Cp|x-y|^2}{(p-1)t}\bigg],\ \ t>0, p>1, f\in \B_b^+(\H).$$
By \cite[Theorem 1.4.1]{Wbook}, this implies that $P_t^0$ has a strictly positive density with respect to $\mu_0$. Therefore, by Theorem \ref{T0},
$\nu\in \scr P_Z$ is unique, and the density $\rr$ has a strictly positive version. In particular, the marginal density $\rr_0$ has a strictly positive version as well. This together with the   Poincar\'e inequality \eqref{P} implies $\log \rr_0\in L^2(\mu_0)$, see step (a) in the proof of \cite[Theorem 2.3(1)]{W16}.

In conclusion, we have the following result, where $\ll_{\kk,\tau}$ and $q_\ll$ are given before Theorem \ref{T3.1}.

\beg{thm}\label{T3.9} In the above framework, let $\kk=\ff 1{\ll_1}$. If $\mu(\e^{\ll |Z|^2})<\infty$ for some $\ll>\ll_{\kk,\tau}$, the $\scr P_Z$ contains a unique measure $\nu=\rr\mu$, where $\rr$ is strictly positive such that $\eqref{SI4}$ and $\eqref{SI3}$ hold for any $p\in (1,\ff{q_\ll}{q_\ll-1})$ and some constant $k=k(\ll,p)>0$.   \end{thm}

\paragraph{Acknowledgement.} The author would like to thank Professor Shige Peng for valuable conversations and the referees for helpful comments and a number of corrections.

\end{document}